\documentclass[12pt]{amsart}
\usepackage[mathscr]{eucal} 
\usepackage{amsmath,amsfonts} 
\parskip=\smallskipamount 
\newtheorem{theorem}{Theorem}[section] 
\newtheorem{proposition}[theorem]{Proposition} 
\newtheorem{corollary}[theorem]{Corollary} 

\theoremstyle{definition} 
 
\newtheorem{example}[theorem]{Example}

\makeatletter 
\@addtoreset{equation}{section} 
\makeatother

\newcommand{\BB}{{\mathbb B}}
\newcommand{\CC}{{\mathbb C}} 
\newcommand{\NN}{{\mathbb N}} 
 
\newcommand{\ZZ}{{\mathbb Z}} 
\newcommand{\DD}{{\mathbb D}} 
\newcommand{\RR}{{\mathbb R}} 
\newcommand{\TT}{{\mathbb T}}
 
\newcommand{\cA}{{\mathcal A}} 
\newcommand{\cB}{{\mathcal B}} 
 
\newcommand{\cD}{{\mathcal D}} 
\newcommand{\cE}{{\mathcal E}} 
\newcommand{\cF}{{\mathcal F}} 
\newcommand{\cG}{{\mathcal G}} 
\newcommand{\cH}{{\mathcal H}} 
 
\newcommand{\cK}{{\mathcal K}} 
\newcommand{\cL}{{\mathcal L}} 
 
\newcommand{\cN}{{\mathcal N}} 
\newcommand{\cO}{{\mathcal O}}
 
\newcommand{\cR}{{\mathcal R}} 
\newcommand{\cS}{{\mathcal S}}

\newcommand{\bah}{\mathbf{H}}
\newcommand{\dom}{\operatorname{Dom}} 
\newcommand{\Kr}{Kre\u\i n} 
\newcommand{\Na}{Na\u\i mark}
\newcommand{\Ra}{\Rightarrow} 
\newcommand{\ran}{\operatorname{Ran}} 
 
\newcommand{\ra}{\rightarrow} 
\newcommand{\ol}{\overline}

\let\phi=\varphi 
\newcommand{\iac}{\mathrm{i}} 
\renewcommand{\ker}{\operatorname{Ker}} 
 
\newcommand{\de}{\mathrm{d}}
\newcommand{\emath}{\mathrm{e}}

\newcommand{\supp}{\operatorname{supp}}
\renewcommand{\Re}{\operatorname{Re}}
\newcommand{\cl}{\operatorname{Clos}}

\begin{document}
\title[Reproducing Kernel Kre\u\i n Spaces]{A Survey on Reproducing Kernel Kre\u\i n Spaces}
 
\author[A. Gheondea]{Aurelian Gheondea} 

\address{Department of Mathematics, Bilkent University, 06800 Bilkent, Ankara, 
Turkey, \emph{and} Institutul de Matematic\u a al Academiei Rom\^ane, C.P.\ 
1-764, 014700 Bucure\c sti, Rom\^ania} 
\email{aurelian@fen.bilkent.edu.tr \textrm{and} A.Gheondea@imar.ro} 

\begin{abstract} This is a survey on reproducing
kernel \Kr\ spaces and their interplay with 
operator valued Hermitian kernels. Existence
and uniqueness properties are carefully reviewed. The approach used in this survey involves the more abstract, but very useful, concept 
of linearisation or Kolmogorov decomposition, as well as the 
underlying concepts of \Kr\ 
space induced by a selfadjoint operator and that of \Kr\ space continuously 
embedded. The operator range feature of reproducing 
kernel spaces is emphasised. 
A careful presentation of Hermitian kernels on complex regions that 
points out a universality property of the Szeg\"o kernels with respect to 
reproducing kernel \Kr\ spaces of holomorphic functions is included.\medskip

\noindent
\textit{Acknowledgement.} Work supported by a grant of the Romanian 
National Authority for Scientific Research, CNCS Ð UEFISCDI, project number 
PN-II-ID-PCE-2011-3-0119.\medskip

\noindent \textit{Key words and phrases.} \Kr\ space, reproducing kernel, Hermitian kernel, linearisation, Kolmogorov decomposition, induced space, holomorphy.\medskip

\noindent \textit{2010 Mathematics Subject Classification.} 46C07, 46C20, 47B32, 47B50, 32A25.

\end{abstract}

\maketitle

\section{Introduction}\label{s:i}

A kernel $K$ on a set $X$ and with operator entries is the analog
of an abstract matrix double indexed on $X$ and for which each entry is a bounded
linear operator acting between appropriate spaces. The adjoint $K^\sharp$ 
of the kernel $K$ can be defined by analogy with the adjoint of matrices. 
A Hermitian kernel $K$ is defined by the property
$K=K^\sharp$. Roughly speaking, a reproducing kernel \Kr\ space on a set $X$ is a 
\Kr\ space $\cR$ of functions on $X$ for which there exists a Hermitian kernel with the 
property that the evaluations of the functions in $\cR$ can be calculated in terms of 
the kernel $K$. One of the main problems is that of 
associating a reproducing kernel \Kr\ space $\cR$ to a given Hermitian kernel $K$.

The classical theory says that 
for any reproducing kernel Hilbert space $\cR$ its reproducing kernel $K$
is positive semidefinite and uniquely determined by $\cR$. Conversely, to any
positive semidefinite kernel $K$ there is associated a unique reproducing kernel 
Hilbert space such that $K$ is its reproducing kernel. These facts can be proven 
also for a slightly more general situation relating in a similar fashion reproducing 
kernel Pontryagin spaces to Hermitian kernels with finite negative (or
positive) signatures.

For genuine indefinite Hermitian kernels the correspondence with reproducing kernel 
\Kr\ spaces is much more complicated. There is a variety of characterisations of 
those Hermitian kernels that produce reproducing kernel \Kr\ spaces in terms of 
boundedness with respect to positive semidefinite kernels and in terms of 
decomposability as a difference of two positive semidefinite kernels. Uniqueness of 
the reproducing kernel \Kr\ space, provided that it exists, is also problematic. 
Characterisations of uniqueness are available in terms of lateral spectral gaps as 
well as in terms of maximal uniformly positive subspaces.

When Hermitian kernels on regions in either a single complex variable or 
several complex 
variables have certain holomorphy properties, existence of reproducing kernel \Kr\ 
spaces is guaranteed. This fact is related to a certain universality property of the 
Szeg\"o kernel that allows the construction of the corresponding reproducing 
kernel \Kr\ space inside the Hardy space or, respectively, the Drury-Arveson space.

The approach used in this survey involves a more abstract, but very useful, concept 
of linearisation or Kolmogorov decomposition, the underlying concept of \Kr\ 
spaces induced by selfadjoint operators, as well as the concept of \Kr\ spaces 
continuously embedded, pointing out the operator range feature of reproducing
kernel spaces.

\section{\Kr\ Spaces and their Linear Operators}
A \emph{\Kr\ space} $\mathcal{K}$ is a complex linear space on 
which it is defined an indefinite scalar product $[\cdot,\cdot]$ such that
$\mathcal{K}$ is decomposed in a direct sum 
\begin{equation} \label{mak}
\mathcal{K} = \mathcal{K}_{+} [\dot{+}] \mathcal{K}_{-} 
\end{equation} 
in such a way that $\mathcal{K}_{\pm}$ are Hilbert spaces with 
scalar products $\pm [\cdot,\cdot]$, and the direct sum in 
\eqref{mak} is orthogonal with respect to the indefinite scalar product 
$[\cdot,\cdot]$, i.e.\ $\mathcal{K}_{+} \cap \mathcal{K}_{-} = \{0\}$ and 
$[x_{+}, x_{-}] = 0$ for all $x_{\pm} \in \mathcal{K}_{\pm}$. 
The decomposition \eqref{mak} gives rise to a positive definite scalar 
product $\langle \cdot,\cdot\rangle$ by setting
$\langle x,y\rangle   
: = \langle x_{+},y_{+}\rangle  - \langle x_{-},y_{-}\rangle$, where $x =
x_{+} + x_{-}, y =  
y_{+} + y_{-}$, and $x_{\pm}, y_{\pm} \in \mathcal{K}_{\pm}$. The 
scalar product $\langle \cdot ,\cdot\rangle $ 
defines on $\mathcal{K}$ a structure of 
Hilbert space. Subspaces $\mathcal{K}_{\pm}$ are orthogonal with 
respect to the scalar product 
$\langle \cdot,\cdot\rangle $, too. One denotes
by $P_{\pm}$ the corresponding orthogonal projections onto 
$\mathcal{K}_{\pm}$, and let $J = P_{+} - P_{-}$. The operator $J$ 
is a \emph{symmetry}, i.e.\ a selfadjoint and unitary operator,
$J^{*} J = J J^{*} = J^{2} = I$. Unless otherwise specified, any \Kr\ space is
considered as a Banach space with the norm given by an arbitrary norm induced
by the positive definite inner product associated to a fundamental
decomposition. For two \Kr\ spaces $\cK_1$
and $\cK_2$, we denote by $\cL(\cK_1,\cK_2)$ the Banach
space of all bounded linear operators $T\colon \cK_1\ra\cK_2$.   

Given a \Kr\ space $(\cK,[\cdot,\cdot])$, the cardinal numbers
\begin{equation}\label{granks} \kappa^+(\cK)=\dim(\cK^+),\quad
  \kappa^-(\cK)=\dim(\cK^-),\end{equation}
do not depend on the fundamental decomposition and they are called,
respectively, the \emph{geometric ranks of positivity/negativity} of
$\cK$.

The operator $J$ is called a \emph{fundamental symmetry} of the 
\Kr\ space $\mathcal{K}$. Note that $[x,y] = \langle J x,y\rangle$,  
$(x,y \in \mathcal{K})$. If $T$ is a densely defined operator from a \Kr\
space $\cK_{1}$ to another \Kr\ space $\mathcal{K}_{2}$, 
\emph{the adjoint} of $T$ can be
defined as an operator $T^{\sharp}$ defined on the set of all 
$y \in \mathcal{K}_{2}$ for 
which there exists $h_{y} \in \mathcal{K}_{1}$ such that $[T 
x,y] = [x, h_{y}]$, and $T^{\sharp} y = h_{y}$. In addition,
$T^{\sharp} = J_{1} T^{*}J_{2}$, where $T^{*}$ denotes the 
adjoint operator of $T$ with respect to the Hilbert spaces
$(\mathcal{K}_1,\langle\cdot,\cdot\rangle_{J_1})$ and
$(\mathcal{K}_2,\langle\cdot,\cdot\rangle_{J_2})$. The symbol $\sharp$ denotes
the adjoint when at least of the spaces $\cK_1$ or $\cK_2$ is
indefinite. 
In the case of an operator $T$ defined on 
the \Kr\ space $\mathcal{K}$, $T$ is called \emph{symmetric} 
if $T \subseteq T^{\sharp},$ i.e.\ if the relation 
$[T x,y] = [x, T y]$ holds for each $x, y \in \dom(T)$, and 
$T$ is called selfadjoint if $T =T^{\sharp}$. 

In this survey, a bit of  geometry of \Kr\ spaces is used. Thus,
a (closed) subspace $\cL$ of a \Kr\ space $\cK$ is called \emph{regular} if
$\cK=\cL+\cL^\perp$, where $\cL^\perp=\{x\in\cK\mid [x,y]=0\mbox{ for all
}y\in\cL\}$. Regular spaces of \Kr\ spaces are important since they are
exactly the analogue of \Kr\ subspaces, that is, if we want $\cL$ be a \Kr\
space with the restricted indefinite inner product and the same strong
topology, then it should be regular. 

In addition, let us recall that, given a subspace $\cL$ of a \Kr\ space, one
calls $\cL$ \emph{non-negative} (\emph{positive})
if the inequality $[x,x]\geq 0$ holds for $x\in \cL$ 
(respectively, $[x,x]>0$ for all $x\in\cL\setminus\{0\}$). 
Similarly one defines \emph{non-positive} and
\emph{negative} subspaces. A subspace $\cL$ is called \emph{degenerate} if
$\cL\cap\cL^\perp\neq \{0\}$. Regular subspaces are non-degenerate. 
As a consequence of the Schwarz inequality, 
if a subspace $\cL$ is either positive or negative it is nondegenerate. A
remarkable class of subspaces are those regular spaces that are either positive
or negative, for which the terms \emph{uniformly positive},
respectively, \emph{uniformly negative} are used. These notions can be defined
for linear manifolds also, that is, without assuming closedness.

A linear operator $V$ defined from a subspace of a \Kr\ space $\cK_1$ and
valued into another \Kr\ space $\cK_2$ is called \emph{isometric} if
$[Vx,Vy]=[x,y]$ for all $x,y$ in the domain of $V$. Note that
isometric operators between genuine \Kr\ spaces may be unbounded and different
criteria of boundedness are available. 
However, in this presentation, a \emph{unitary} operator between \Kr\
spaces means that it is a bounded isometric operator that have bounded inverse.

\section{Hermitian Kernels}\label{s:hk}
Let $ X $ be a nonempty set and let $\bah=\{ \cH_x\}_{x\in X }$ be a
family of \Kr\ spaces with inner products denoted by
$[\cdot,\cdot]_{ \cH_x}$. A mapping $K$ defined on $ X \times X $
such that $K(x,y)\in \cL( \cH_y, \cH_x)$ for all $x,y\in X $ is
called an \emph{$\bah$-kernel on $X$}. In case $\cH_x=\cH$ for all $x\in X$, 
where $\cH$ is some fixed \Kr\ space, one talks about an $\cH$-kernel on $X$, 
while, even more particularly, if $\cH=\CC$, one talks about 
a \emph{scalar kernel on $X$}, or simply a \emph{kernel} on $X$.

To any $\bah$-kernel $K$ one associates its 
\emph{adjoint} $K^\sharp$ defined by $K^\sharp(x,y)=K(y,x)^\sharp$, 
for all $x,y \in X$. The 
$\bah$-kernel\ $K$ is called \emph{Hermitian} if $K=K^\sharp$, that is,
 \begin{equation}\label{e:anul} K(x,y)=K(y,x)^\sharp,\quad x,y\in X . \end{equation}

Denote by $ \cF(\bah)$ the set of all $\bah$-vector fields $f$ on $X$, that is,
$f=\{f_x\}_{x\in X }$ such that $f_x\in \cH_x$, for all
$x\in X $, and let $\cF_0(\bah)$ denote the set of all $f\in \cF(\bah)$
of finite support, that is, the set $\supp(f)=\{x\in X \mid f_x\neq0\}$ is finite. 
Alternatively, one can view any $f\in\cF(\bah)$ as a function 
$f\colon X\ra \bigcup_{x\in X}\cH_x$ such that $f(x)\in \cH_x$ for all $x\in X$.

If $K$ is a Hermitian $\bah$-kernel then one can introduce on $ \cF_0(\bah)$ an
inner product $[\cdot,\cdot]_K$ defined by
 \begin{equation}\label{e:fega} [f,g]_K=\sum_{x,y\in X }
 [ K(x,y)f(y),g(x)]_{ \cH_x},\quad f,g\in \cF_0(\bah). \end{equation}
The $\bah$-kernel $K$ is called \emph{positive semidefinite} if
 \begin{equation} \sum_{x,y\in X }[ K(x,y)h(y),h(x)]_{ \cH_x}\geq0,
\quad h\in \cF_0(\bah). \end{equation}
Every positive
semidefinite  $\bah$-kernel\  is Hermitian. Also, a
Hermitian  $\bah$-kernel\ is positive semidefinite if and only if the 
corresponding inner product in (\ref{e:fega}) is nonnegative.

Let us denote by  ${\mathfrak K}^h(\bah)$ the class of all Hermitian
$\bah$-kernels
and by ${\mathfrak K}^+(\bah)$ the subclass of all positive semidefinite
$\bah$-kernels. On
${\mathfrak K}^h(\bah)$ one defines addition, subtraction and multiplication
with real numbers in a natural way. Moreover, on ${\mathfrak K}^h( \bah)$ one
has a natural partial order defined as follows: if $H,K\in{\mathfrak
K}^h(\bah)$ then $H\leq K$ means $[f,f]_H\leq [f,f]_K$, for all
$f\in \cF_0(\bah)$. With this definition one gets
 \begin{equation}\label{e:kaplus} 
 {\mathfrak K}^+(\bah)=\{H\in{\mathfrak K}^h(\bah)\mid H\geq 0\}, 
 \end{equation} and
${\mathfrak K}^+(\bah)$ is a strict cone of ${\mathfrak K}^h(\bah)$, that is, it
is closed under addition and multiplication with nonnegative numbers,
and ${\mathfrak K}^+(\bah)\bigcap -{\mathfrak K}^+(\bah)=\{0\}$, where $0$ denotes 
the null kernel.

More generally, one can define the 
\emph{signatures} of $K$, denoted by
$\kappa_\pm(K)$, as the positive/negative signatures of the inner product 
$[\cdot,\cdot]_K$. Then, if $\kappa(K)=\min\{\kappa_-(K),\kappa_+(K)\}$ denotes
the \emph{definiteness signature} of $K$, the Hermtian $\cH$-kernel $K$ is called 
\emph{quasi semidefinite} if $\kappa(K)<\infty$, that is, 
either $\kappa_-(K)$ or $\kappa_+(K)$ is finite.

A pairing $(\cdot,\cdot)_\cF$ can be defined for
arbitrary $f,g\in \cF(\bah)$, provided at least one of $f$ and $g$ has
finite support, by
\begin{equation}\label{e:inop}(f,g)_\cF
=\sum_{x\in X }[f(x),g(x)]_{ \cH_x}. \end{equation}
When restricted to $\cF_0(\bah)$ this pairing becomes a nondegenerate 
inner product.
To each  $\bah$-kernel $K$ one associates the {\em convolution operator},
denoted also by $K$, and defined by
\begin{equation}\label{e:convol}K\colon \cF_0(\bah)\ra \cF(\bah),\quad
(Kf)(x)=\sum_{y\in X }K(x,y)f(y),\quad f\in \cF_0(\bah). \end{equation}
Then 
\begin{equation}\label{e:conv} [f,g]_K=(Kf,g)_\cF,\quad f,g\in\cF_0(\bah).
\end{equation}
Consequently, the kernel $K$ is positive semidefinite (Hermitian) if and only if the
corresponding convolution operator $K$ is positive semidefinite (Hermitian), that is,
$(Kf,g)_\cF\geq 0$ ($(Kf,g)_\cF=(f,Kg)_\cF$), 
for all $f,g\in \cF_0(\bah)$. Similar assertions can
be made about the signatures $\kappa_\pm(K)$, with appropriate 
definitions of signatures of Hermitian operators on inner product spaces.

At this point, it is worth noting that 
there is no restriction of generality if one 
assumes that all spaces $\cH_x$ are Hilbert. 
To see this, fixing a fundamental symmetry $J_x$ on each \Kr\ space $\cH$, one
can refer to the Hilbert spaces $(\cH_x;\langle\cdot,\cdot\rangle_{J_x})$ and, 
one considers the kernel $H$ on $X$ defined by $H(x,y)=J_x K(x,y)J_y$, for all 
$x,y\in X$. Taking into account that $A^\sharp=J_x A^* J_y$ for any bounded 
linear 
operator $A\colon \cH_x\ra\cH_y$, it follows that all notions like Hermitian, 
positive semidefinite, signatures, etc.\ defined for the kernel $K$, 
have word-for-word
transcriptions for the kernel $H$. However, since this survey is focusing on the 
indefinite case, allowing right from the beginning \Kr\ spaces $\cH_x$ brings more 
symmetry and simpler formulas.

A \emph{reproducing kernel inner product space}, with respect to the set $X$, the 
collection of \Kr\ spaces $\bah=\{\cH_x\}_{x\in X}$, 
and $\bah$-kernel $K$, is an inner product space 
$(\cR,[\cdot,\cdot]_\cR)$ subject to the following conditions:

\begin{itemize}
\item[(rk1)] $\cR\subseteq\cF(\bah)$.
\item[(rk2)] For all $y\in X$ and all $h\in \cH_y$ the map 
$X\ni x\mapsto K(x,y)h\in  \cH_x$ belongs to $\cR$.
\item[(rk3)] $[ f(x),h]_{\cH_x}=[f,K(\cdot,x)h]_{\cR}$, for all $f\in\cR$,
$x\in X$, and $h\in\cH_x$.
\end{itemize}
The axiom (rk3) is usually called the \emph{reproducing property} while the 
$\bah$-kernel $K$ is called the \emph{reproducing kernel} of $\cR$. Note
that $K$ is necessarily a Hermitian $\bah$-kernel.

According to the axiom (rk2), it is useful to consider the notation 
$K_y=K(\cdot,y)\colon X\ra \cL(\cH_y,\cR)$, for $y\in X$. 
An immediate consequence of the axioms (rk1)--(rk3) is that the inner product 
$[\cdot,\cdot]_\cR$ is nondegenerate. Also, recalling that on any 
nondegenerate inner product space the weak topology is separated, 
the following 
\emph{minimality property} holds
\begin{itemize}
\item[(rk4)] The span of $\{K_x\cH_x\mid x\in X\}$ is weakly dense in $\cR$.
\end{itemize}
Also, the reproducing kernel is uniquely determined by the reproducing 
kernel inner product space.

Given a Hermitian $\bah$-kernel $K$ on $X$, consider 
the subspace $\cR_0(\bah)$ of $\cF(\bah)$ spanned by $K_xh$, 
for all $x\in X$ and all $h\in \cH_x$, on which one can define the inner product
\begin{equation}\label{e:hega} 
[\sum_{i=1}^m K_{x_i} h_i,\sum_{j=1}^nK_{y_j} k_j]_{\cR_0}
=\sum_{i=1}^m\sum_{j=1}^n[K(y_j,x_i)h_i,k_j]_{\cH_{y_j}},
\end{equation}
This definition can be proven to be correct:  vectors in $\cR_0(\bah)$ may have 
different representations as linear combinations of $K_yk$ but the definition in 
\eqref{e:hega} is independent of these. The subspace $\cR_0(\bah)$ of $\cF(\bah)$ 
is the range of the convolution operator $K$ defined at \eqref{e:convol} and
the inner product $[\cdot,\cdot]_{\cR_0}$ is nondegenerate. In addition, 
$(\cR_0(\bah);[\cdot,\cdot]_{\cR_0})$ is a 
reproducing kernel inner product space with reproducing kernel $K$. 

In case the reproducing kernel inner product space
$(\cR;[\cdot,\cdot]_\cR)$ is a \Kr\ space, one talks about a \emph{reproducing 
kernel \Kr\ space}. The uniqueness of the reproducing kernel of a \Kr\ space has
a stronger characterisation.

\begin{theorem}\label{t:reka} Let $\cK$ be a \Kr\ space of\, $\bah$-valued vector 
fields on $X$, that is, $\cK\subseteq\cF(\bah)$. For each $x\in X$ 
consider the linear operator $E(x)\colon \cK\ra\cH_x$ of evaluation at $x$, that is, 
$E(x)f=f(x)$ for all $f\in\cK$. 
Then, $\cK$ has a reproducing kernel if and only if $E(x)$ is bounded for all 
$x\in X$. 
In this case, the reproducing kernel of $\cK$ is
\begin{equation}\label{e:heka} K(x,y)=E(x)E(y)^\sharp,\quad x,y\in X,
\end{equation} hence uniquely determined by $\cK$.
\end{theorem}
With the notation as in Theorem~\ref{t:reka}, it is useful to note that, if $\cK$ is 
a reproducing kernel \Kr\ space with reproducing kernel $K$, then $K_x$ can be 
viewed as a linear operator $\cH_y\ra \cK$ and $K_x=E(x)^\sharp$ for all $x\in X$. 
In particular, 
\begin{equation}\label{e:kaka} K(x,y)=K_x^\sharp K_y,\quad x,y\in X.
\end{equation}

Conversely, 
given a Hermitian $\bah$-kernel $K$ on $X$, the questions on existence and 
uniqueness of a reproducing kernel \Kr\ space $\cR\subseteq\cF(\bah)$ such 
that $K$ is its reproducing kernel are much more difficult. 
If $(\cR;[\cdot,\cdot]_\cR)$ is a reproducing kernel \Kr\ space with 
reproducing kernel $K$,
in view of the axioms (rk1)--(rk3) and the minimality property (rk4), the inner 
products $[\cdot,\cdot]_\cR$ and $[\cdot,\cdot]_{\cR_0}$ coincide on $\cR_0(\bah)$ 
and $\cR_0(\bah)$ is dense in $\cR$. Thus, existence of a reproducing kernel 
\Kr\ space associated to a given $\bah$-kernel $K$ on $X$ depends heavily on the 
possibility of "completing" the inner product space 
$(\cR_0(\bah);[\cdot,\cdot]_{\cR_0})$ to a \Kr\ space $\cR$ inside $\cF(\bah)$, which 
is a core problem in the theory of indefinite inner product spaces.

In order to tackle these questions, let us note that if $\cR$ is a reproducing kernel 
Hilbert space with reproducing kernel $K$, then $K$ should be positive semidefinite.
Conversely, if $K$ is a positive semidefinite $\bah$-kernel on $X$, then the inner 
product space $(\cR_0,[\cdot,\cdot]_{\cR_0})$ defined at \eqref{e:hega} is a 
pre-Hilbert space and the existence of a reproducing kernel Hilbert space $\cR$ 
with reproducing kernel $K$ depends on whether $(\cR_0,[\cdot,\cdot]_{\cR_0})$
has a completion inside of $\cF(\bah)$. 
Actually, this is always the case and, 
moreover, uniqueness holds as well, but these two facts are slightly more 
general, 
namely, they are true if $K$ is quasi semidefinite, in which case it is
obtained a  
unique reproducing kernel Pontryagin space, cf.\ Section~\ref{s:qsk}.

In this survey, existence and uniqueness of reproducing kernel \Kr\ spaces 
associated to Hermitian kernels are approached through the more abstract, 
but very useful, concept of linearisation or Kolmogorov decomposition.
By definition, a \emph{linearisation}, sometimes called 
\emph{Kolmogorov decomposition}, of the $\bah$-kernel $K$ is a
pair $(\cK;V)$, subject to the following conditions:
\begin{itemize}
\item[(kd1)] $\cK$ is a   \Kr\ space and $V(x)\in\cL(\cH_x,\cK)$ for all $x\in X$.
\item[(kd2)] $K(x,y)=V(x)^\sharp V(y)$ for all $x,y\in X$.
\end{itemize}
The linearisation $(\cK;V)$ is called \emph{minimal} if the following condition holds 
as well:
\begin{itemize}
\item[(kd3)] $\cK=\bigvee\limits_{x\in X} V(x)\cH_x$.
\end{itemize}

Two linearisations $(\cK_j;V_j)$, $j=1,2$, of the same $\bah$-kernel $K$ are 
called \emph{unitary equivalent} if there exists a   \Kr\ space bounded unitary 
operator $U\colon \cK_1\ra \cK_2$ such that $V_2(x)=UV_1(x)$ for all $x\in X$.

There is a close connection between the notion of reproducing kernel  \Kr\ 
space with reproducing $\bah$-kernel $K$ and that of minimal linearisation of $K$.

\begin{proposition}\label{p:lik} \emph{(1)} 
Let $\cR$ be a reproducing kernel   \Kr\ space with reproducing 
$\bah$-kernel $K$. Then, letting 
\begin{equation}\label{e:vex}V(x)=K_x=K(\cdot,x)\in\cL(\cH_x,\cR),
\end{equation} the pair
$(\cR;V)$ is a minimal linearisation of the $\bah$-kernel $K$.

Conversely, letting $(\cK;V)$ be a minimal linearisation of the $\bah$-kernel $K$, 
then,
\begin{equation}\label{e:res} \cR=\{V(\cdot)^\sharp\mid k\in\cK\}, 
\end{equation}
is a vector subspace of $\cF(\bah)$ which, with respect to the
the inner product defined by
\begin{equation}\label{e:ipres} [V(\cdot)^\sharp h,V(\cdot)^\sharp k]_\cR=[h,k]_\cK,\quad 
h,k\in \cK,
\end{equation} is a \Kr\ space with reproducing kernel $K$.

\emph{(2)} In the correspondence defined at \eqref{e:res} and \eqref{e:ipres}, 
two unitary 
equivalent minimal realisations of the same $\bah$-kernel $K$ produce the same 
reproducing kernel Krein space and hence, the correspondence between 
reproducing kernel Krein spaces and minimal linearisations is one-to-one, 
provided that unitary equivalent minimal linearisations are identified.
\end{proposition}

\section{Some Examples}\label{s:se}

\begin{example} \label{ex:matrices}
\emph{Matrices.} Let $X=\{1,2,\ldots,n\}$ and $\cH_x=\CC$ 
for each $x\in X$. A kernel on $X$ is simply a map $K\colon X\times X\ra \CC$ and 
hence, it can be identified with the $n\times n$ complex matrix $[K(i,j)]_{i,j=1}^n$. 
The 
kernel $K$ is Hermitian if and only the matrix $[K(i,j)]_{i,j=1}^n$ is Hermitian. The
vector space $\cF=\cF_0=\CC^n$ hence, in view of the definition \eqref{e:convol},
the convolution operator $K\colon \cF\ra\cF$ is simply the linear operator 
$K\colon \CC^n\ra\CC^n$ associated to the matrix $K$. 

If $K$ is a Hermitian kernel on $X$ then the inner product $[\cdot,\cdot]_K$ is
the familiar inner product associated to the Hermitian matrix $[K(i,j)]_{i,j=1}^n$. 
The signatures $\kappa_\pm(K)$ coincide, respectively, with the number of 
positive/negative eigenvalues of the matrix $[K(i,j)]_{i,j=1}^n$, counted with 
multiplicities. $K$ is positive semidefinite if and only if the matrix $[K(i,j)]_{i,j=1}^n$ is 
positive semidefinite.

Assuming that the kernel $K$ on $X$ is Hermitian, for each $j=1,\ldots,n$, the "map" 
$K_j\colon X\ra \CC$ is simply the column vector $[K(i,j)]_{i=1}^n$ in $\CC^n=\cF$.
The vector space $\cR_0$, the range of the convolution operator $K$, is the vector
subspace of $\CC^n$ generated by all the column vectors $[K(i,j)]_{i=1}^n$. In this 
particular case, $\cR_0$ is the reproducing kernel \Kr\ space with reproducing 
kernel $K$. The inner product $[\cdot,\cdot]_{\cR_0}$ is defined at \eqref{e:hega}.
 \end{example}

\begin{example}\label{ex:opm} 
\emph{Operator Block Matrices.} Let $X=\{1,2,\ldots,n\}$ but, this
time, for each $x\in X$ one denotes by $\cH_x$ a Hilbert space and let
$\bah=\{\cH_x\}_{x\in X}$. An $\bah$-kernel on $X$, originally defined as a map $K$
on $X\times X$ such that $K(x,y)\in \cL(\cH_y,\cH_x)$, is naturally identified with the
operator block matrix $[K(i,j)]_{i,j=1}^n$. The vector space $\cF(\bah)=\cF_0(\bah)$ 
together with its natural inner product $(\cdot,\cdot)_\cF$ defined at \eqref{e:inop}
is naturally identified with the Hilbert space $\cH=\bigoplus_{i=1}^n \cH_i$. Then, 
the convolution operator $K$ associated to the $\bah$-kernel $K$ is identified with
the bounded
linear operator $K\colon \cH\ra \cH$ naturally associated to the operator block
matrix $[K(i,j)]_{i,j=1}^n$. For each $j=1,\ldots,n$, the map $K_j=K(\cdot,j)$ is
the operator block column matrix $[K(i,j)]_{i=1}^n\colon \cH_j\ra \cH$.

The $\bah$-kernel $K$ on $X$ is Hermitian if and only if the corresponding operator
block matrix $[K(i,j)]_{i,j=1}^n$ is Hermitian, if and only if the convolution operator 
$K\in \cL(\cH)$ is selfadjoint (Hermitian). In this case, the vector space 
$\cR_0(\bah)$ is identified with the range of the convolution operator $K$ as a 
subspace of $\cH$ and 
is spanned by the ranges of the operator block column matrices 
$K_j=[K(i,j)]_{i=1}^n\colon \cH_j\ra \cH$, for $j=1,\dots,n$. 
The inner product $[\cdot,\cdot]_{\cR_0}$ is defined as in \eqref{e:hega}.
\end{example}

\begin{example}\label{ex:hardy}
\emph{The Hardy space $H^2(\DD)$.} Let $\DD$ denote the 
open unit ball in the complex field and consider the \emph{Szeg\"o kernel}
\begin{equation}S(z,w)=\frac{1}{1-z\overline{w}}=\sum_{n=0}^\infty \ol{w}^nz^n,\quad 
z,w\in\DD,
\end{equation} the series converging absolutely and 
uniformly on any compact subset of 
$\DD\times\DD$.
Clearly, $S$ is a Hermitian scalar kernel on $\DD$ and for every $w\in \DD$
\begin{equation*} S_w(z)=S(z,w)=\sum_{n=0}^\infty \ol{w}^nz^n,\quad 
z,w\in\DD,
\end{equation*} the series converging absolutely and uniformly on any compact 
subset of $\DD$, $S_w$ is a holomorphic function on $\DD$, and its Taylor
coefficients are $1,\ol{w},\ol{w}^2,\ldots$. In order to describe the reproducing kernel 
inner product space associated to the Szeg\"o kernel $S$, in view of \eqref{e:hega},
one has to consider the vector space generated by the functions $S_w$, 
$w\in \DD$, and complete it with respect to the inner product 
\begin{equation}\label{e:suw}\langle S_u,S_w\rangle= S_u(w)=S(w,u)
=\sum_{n=0}^\infty \ol{u}^nw^n,
\quad u,v\in\DD.\end{equation}
 This completion is called the \emph{Hardy space}, denoted by $H^2(\DD)$, and in 
 view of \eqref{e:suw}, is
 \begin{equation}\label{e:hardy} H^2(\DD)=\{f\mid f(z)=\sum_{n=0}^\infty a_n z^n,\ 
 \sum_{n=0}^2 |a_n|^2<\infty\}.
 \end{equation} Since $H^2(\DD)$ is a Hilbert space, it follows that the Szeg\"o 
 kernel is positive semidefinite.
 
The Hardy space $H^2(\DD)$ has some other special properties, among which, a 
distinguished role is played by the boundary values of its functions. More precisely, 
for each $f\in H^2(\DD)$, there exists $\widetilde f$ a function defined on the unit 
circle $\TT=\partial \DD$, such that
\begin{equation}\widetilde f(\emath^{\iac t})=\lim_{r\ra 1-}f(r\emath^{\iac t}),\quad 
\mbox{ a.e. } t\in[0,2\pi).
\end{equation} The function $\widetilde f$ is uniquely determined by $f$, a.e.\ on 
$\TT$, and $\widetilde f\in L^2(\TT)$. Usually, there is no distinction between a 
function $f$ in 
$H^2(\DD)$ and $\widetilde f$, that is, the tilde sign is not used at all.
In this way, an isometric embedding of 
$H^2(\DD)$ into $L^2(\TT)$ is defined, in particular,
\begin{equation}
\label{e:hep} \langle f,g\rangle_{H^2(\DD)}=\langle f,
g\rangle_{L^2(\TT)}=\frac{1}{2\pi} \int_0^{2\pi}  f(\emath^{\iac t})\overline{
g(\emath^{\iac t})}\de t,\quad f,g\in H^2(\DD).\end{equation}
With respect to this embedding, $H^2(\DD)$ is identified with the 
subspace of $L^2(\TT)$ of all functions $f$ whose Fourier coefficients of 
negative index vanish. 
\end{example}

\begin{example}\label{ex:das}
\emph{The Drury-Arveson Space.}
Let $\BB_r(\xi)$ be the open
ball of radius $r$ and center $\xi $ in the Hilbert space $\cG=\CC^N$, with inner 
product $\langle \xi,\eta\rangle_\cG=\ol{\eta}\cdot\xi=
\sum_{n=1}^N \overline\eta_n\xi_n$ for 
$\xi=(\xi_n)_{n=1}^N$ and $\eta=(\eta_n)_{n=1}^N$ in $\CC^N$. We write
$\BB_r$ instead of $\BB_r(0)$. The \emph{Szeg\"o kernel} is
\begin{equation}\label{szego}
S(\xi ,\eta )=\frac{1}{1-\langle \xi,\eta\rangle_\cG},\quad \xi,\eta\in \BB_1,
\end{equation} and note that $S$ is a scalar Hermitian kernel on $\BB_1$.
 We now describe a minimal linearisation of $S$. Let
$F(\cG )=\bigoplus\limits _{n=0}^{\infty }{\cG} ^{\otimes n}$
be the \emph{Fock space} associated to $\cG=\CC^N$, 
where $\cG^{\otimes 0}=\CC$ and $\cG ^{\otimes n}$ is the $n$-fold 
Hilbert space tensor product of $\cG$ with itself, hence $F(\cG)$ is a Hilbert space. 
Let 
\begin{equation}\label{e:peme}P_n=\frac{1}{n!}\sum _{\pi \in S_n}\hat\pi
\end{equation}
be the orthogonal projection of $\cG ^{\otimes n}$ onto its
symmetric part, where 
$\hat\pi (\xi _1\otimes \ldots \otimes \xi _n)=\xi _{\pi ^{-1}(1)}
\otimes \ldots \otimes \xi _{\pi ^{-1}(n)}$ for any element $\pi$ 
of the permutation group $S_n$ on $n$ symbols. 
Recall that the \emph{symmetric Fock space} is 
$F^s(\cG )=(\bigoplus\limits _{n=0}^{\infty }P_n)F(\cG )$.
For $\xi \in \BB_1$ set $\xi ^{\otimes 0}=1$ and let $\xi ^{\otimes n}$ 
denote the $n$-fold tensor product $\xi\otimes \ldots 
\otimes \xi $, $n\geq 1$. Note that
\begin{equation*}\|\bigoplus _{n\geq 0}\xi ^{\otimes n}\|^2_{F(\cG)}=
\sum _{n\geq 0}\|\xi  ^{\otimes n}\|^2_{\cG^{\otimes n}}
=\sum _{n\geq 0}\|\xi \|^{2n}_{\cG}=
\frac{1}{1-\|\xi \|^2_\cG}.\end{equation*}
Hence $\bigoplus\limits _{n\geq 0}\xi ^{\otimes n}\in F^s(\cG )$
and one can define the mapping $V_S$ from $\BB_1$ into  $F^s(\cG )$, 
\begin{equation}V_S(\xi )=\bigoplus _{n\geq 0}\xi ^{\otimes n},\quad
  \xi\in\cG.\end{equation} 

The pair $(F^s(\cG );V_S)$ is a minimal linearisation of the kernel $S$.
In order to see this,
$V_S(\xi )$ is also viewed as a bounded linear operator from 
$\CC $ into $F^s(\cG )$ by $V_S(\xi )\lambda =\lambda V_S(\xi )$,
$\lambda \in \CC $, so that, for $\xi ,\eta \in \BB_1$, 
\begin{equation*}
V_S(\eta )^*V_S(\xi )=\langle V_S(\xi ),V_S(\eta )\rangle_{F^s(\cG)}
=\sum _{n\geq 0}\langle \xi ^{\otimes n},\eta ^{\otimes n}\rangle_{\cG^{\otimes n} }
 =\sum _{n\geq 0}(\ol\eta \cdot\xi) ^n=
\displaystyle\frac{1}{1-\ol\eta \cdot\xi }=S(\xi ,\eta ).
\end{equation*} In particular, this shows that the Szeg\"o kernel is positive 
semidefinite. 
The set $\{V_S(\xi )\mid \xi \in \BB_1\}$ is total in 
$F^s(\cG )$ since for $n\geq 1$ and $\xi\in\cG$ one has
$\displaystyle\frac{d^n}{dt^n}V(t\xi )|_{t=0}=n!\xi ^{\otimes n}$.

The reproducing kernel Hilbert space
associated to the Szeg\"o kernel $S$, called the \emph{Drury-Arveson space} and
denoted by $H^2(\BB_1)$, is given by the completion of the linear space 
generated by the functions $S_{\eta }=S(\cdot ,\eta )$, $\eta \in \BB_1$, 
with respect to the inner product defined by 
$\langle S_{\eta },S_{\xi }\rangle_{H^2(\BB_1)}=
S(\xi ,\eta )$, see \eqref{e:hega}. We use the multiindex notation: for 
$n=(n_1,\ldots,n_N)\in \NN_0^N$, let $|n|=n_1+\cdots+n_N$, $n!=n_1!\cdots 
n_N!$, and $\xi^n=\xi_1^{n_1}\cdots \xi_N^{n_N}$. Then, a function $f$ 
holomorphic in $\BB_1$, with Taylor series representation 
$f(z)=\sum_{n\in\NN_0^N} a_n \xi^n$, belongs to $H^2(\BB_1)$ if and only if
\begin{equation}\label{e:da} \|f\|^2_{H^2(\BB_1)}=\sum_{n\in\NN_0^N} 
\frac{n!}{|n|!} |a_n|^2<\infty.
\end{equation}
Note that there exists a unitary operator $\Phi $ from 
the Drury-Arveson space $H^2(\BB_1)$ onto $F^s(\cG )$ such that
$\Phi S_{\xi }=V_S(\xi )$, $\xi \in \BB_1$. 

If $N=1$, the Drury-Arveson space coincides with the Hardy space $H^2(\DD)$ 
described at the previous example.
\end{example}

\begin{example}\label{ex:sc}\emph{Holomorphic Kernels.} 
Given two \Kr\ spaces $\cG$ 
and  $\cH$ let $\Omega$ be a subregion of the open unit disc $\DD$ in the complex 
plane and $\Theta\colon \Omega\ra\cL(\cG,\cH)$. One considers the following 
kernels
\begin{equation} K_\Theta(z,w)=\frac{I-\Theta(z)\Theta(w)^\sharp}{1-z\ol w},
\end{equation}
\begin{equation} K_{\widetilde\Theta}(z,w)=\frac{I-\widetilde\Theta(z)\widetilde\Theta(w)^\sharp}{1-z\ol w},
\end{equation}
\begin{equation}D_\Theta(z,w)=\left[\begin{matrix} 
K_\Theta(z,w) & \displaystyle\frac{\Theta(z)-\Theta(\ol w)}{z-\ol w} \\ \displaystyle
\frac{\widetilde\Theta(z)-\widetilde\Theta(\ol w)}{z-\ol w} & K_{\widetilde\Theta}(z,w)
\end{matrix}\right],
\end{equation}
where $\widetilde\Theta(z)=\Theta(\ol z)^\sharp$. These are operator valued 
Hermitian holomorphic kernels that are of interest in connection to \emph{Schur 
classes} and their generalisations. 

The classical Schur class corresponds to the case when $\cG$ and $\cH$ 
are Hilbert spaces and the Schur kernel $K_\Theta$ is positive semidefinite. 
In case $\cG$ and $\cH$ are \Kr\, or even genuine Pontryagin spaces, in order 
to define the Schur class, positive semidefiniteness has to be imposed on all 
$K_\Theta$, $K_{\widetilde\Theta}$ and $D_\Theta$. 
Given $\kappa\in\NN$, the \emph{generalized Schur class} 
corresponds to the requirement
that the negative signatures of each kernel 
$K_\Theta$, $K_{\widetilde\Theta}$ and $D_\Theta$ are $\kappa$. 

When $\cG=\cH$ and for appropriate regions in the complex field, the 
\emph{Carath\'eodory kernel} $C_\Theta$
\begin{equation}\label{e:ck} C_\Theta(z,\zeta)
=\frac{1}{2} \frac{\Theta(z)-\Theta(\zeta)^\sharp}{1-\ol\zeta z},
\end{equation} as well as the \emph{Nevanlinna kernel} $N_\Theta$
\begin{equation}\label{e:nk} N_\Theta(z,\zeta)
=\frac{\Theta(z)-\Theta(\zeta)^\sharp}{z-\ol \zeta},
\end{equation} are of special interest. As in the case of the Schur classes, positive
semidefiniteness of the corresponding kernels define the 
\emph{Carath\'eodory class}
and, respectively, the \emph{Nevanlinna class} 
of holomorphic functions. Generalized 
classes correspond to the case when the appropriate kernels have fixed negative 
signatures.
\end{example}

\begin{example}\label{ex:tk} \emph{T\"oplitz Kernels.} Let $\cH$ be a \Kr\ space
and $X=\ZZ$, the set of integer numbers. An $\cH$-kernel on $\ZZ$,  
$K\colon \ZZ\times \ZZ\ra\cL(\cH)$, is called a \emph{T\"oplitz Kernel} if 
$K(i,j)=T(i-j)$ for some function $T\colon \ZZ\ra\cL(\cH)$, called the 
\emph{symbol} of $K$. 

One considers the set ${\mathfrak T}(\cH)$ of all T\"oplitz
$\cH$-kernels. In the following it is considered
the subclass of T\"oplitz Hermitian $\cH$-kernels ${\mathfrak
T}^h(\cH)$ and the subclass of T\"oplitz positive semidefinite
$\cH$-kernels ${\mathfrak T}^+(\cH)$. One notes that ${\mathfrak
T}^h(\cH)$ is closed under addition, subtraction and (left and right)
multiplication with bounded operators on $\cH$. Also, ${\mathfrak
T}^+(\cH)$ is a strict cone of ${\mathfrak
T}^h(\cH)$.

Considering the complex vector space $\cF_0(\cH)$ of all functions
$h\colon\ZZ\ra\cH$ with finite support and for an arbitrary Hermitian
kernel $H\in{\mathfrak K}^h(\cH)$ one associates the inner product space
$(\cF_0,[\cdot,\cdot]_H)$ as in the previous sections.
On the vector space $\cF_0(\cH)$ one considers two operators, the
{\em forward shift} $S_+$ defined by $(S_+h)(n)=h(n-1)$, for all
$h\in\cF_0(\cH)$ and $n\in\ZZ$, and the {\em backward shift} $S_-$
defined by $(S_-h)(n)=h(n+1)$, for all $h\in\cF_0(\cH)$ and all $n\in\ZZ$.
If $H\in{\mathfrak K}^h(\cH)$ then it is a T\"oplitz
kernel if and only if
\begin{equation*}[S_+h,g]_H=[h,S_-g]_H,\quad f,g\in\cF_0(\cH).
\end{equation*}
If $H$ is a Hermitian T\"oplitz $\cH$-kernel then both $S_+$ and
$S_-$ are isometric with respect to the inner product $[\cdot,\cdot]_H$,
that is, for all $h,g\in\cF_0(\cH)$  the equalities $[S_+h,S_+g]_H=[h,g]_H$
and $[S_-h,S_-g]_H=[h,g]_H$ hold. The converse is also true, if either $S_+$
or $S_-$ is isometric with respect to the Hermitian $\cH$-kernel $H$
then $H$ is T\"oplitz.

Let $H$ be a T\"oplitz Hermitian $\cH$-kernel. A {\em \Na\ dilation} of $H$
is, by definition, a triple $(U,Q;\cK)$ with the following
properties:
\begin{itemize}
\item[(nd1)] $\cK$ is \Kr\ space, $U\in\cL(\cK)$ is a
unitary operator, and $Q\in\cL(\cH,\cK)$.\smallskip
\item[(nd2)] $\cK=\bigvee\limits_{n\in\ZZ}U^nQ\cH$.
\item[(nd3)] $H(i,j)=Q^\sharp U^{i-j}Q$,
$i,j\in\ZZ$.
\end{itemize}
If the T\"oplitz Hermitian $\cH$-kernel $H$ has
a \Na\ dilation  $(U,Q;\cK)$, letting
$V(n)=U^nQ\in\cL(\cH,\cK)$ it is readily verified that the pair
$(V;\cK)$ is a linearisation of $H$.
\end{example}

\section{Quasi Semidefinite Kernels}\label{s:qsk}
In the following, for simplicity and because most of the applications presented in this survey do not require the
full generality as considered in Section~\ref{s:hk}, there is consider only the case of 
$\cH$-kernels, that is, $\cH_x=\cH$ for some fixed \Kr\ space $\cH$. A Hermitian 
kernel $H\colon X\times X\ra \cL(\cH)$ is associated to an indefinite inner product 
space $[\cdot,\cdot]_H$ on $\cF(\cH)$, the vector space of all complex 
valued functions $f\colon X\ra \cH$, as in \eqref{e:fega}. Recall that 
$\kappa_\pm(H)$, the signatures of $H$, are defined as the positive/negative 
signatures of the inner product $[\cdot,\cdot]_H$, while its \emph{rank of 
indefiniteness} is defined by $\kappa(H)=\min\{\kappa_-(H),\kappa_+(H)\}$. 
A Hermitian $\cH$-kernel 
$H$ is \emph{quasi semidefinite} if $\kappa(H)<\infty$, that is, 
either $\kappa_-(H)$ or $\kappa_+(H)$ is finite.

\begin{theorem}\label{t:qsk} Let $H$ be quasi semidefinite $\cH$-kernel on $X$. 
Then:
\begin{itemize}
\item[(a)] $H$ admits a minimal linearisation to a Pontryagin space $(\cK;V)$ with 
$\kappa(H)=\kappa(\cK)$, unique up to a unitary equivalence.
\item[(b)] There exists a unique reproducing kernel Pontryagin space $\cR$ on $X$
and with reproducing kernel $H$.
\end{itemize}
\end{theorem}

On the vector space $\cF(\cH)$ of $\cH$-fields on $X$, the inner product 
$[\cdot,\cdot]_H$, see \eqref{e:fega} is defined. Since $\kappa(H)$ is finite, the
inner product space $(\cF(\cH);[\cdot,\cdot]_H)$ is decomposable. 
To make a choice, assume that $\kappa_-(H)<\infty$, hence
\begin{equation}\label{e:dec} \cF(\cH)=\cF_-[+]\cF_0[+]\cF_+,
\end{equation} where $\cF_0$ is the isotropic subspace, $\cF_\pm$ 
are positive/negative subspaces, and $\dim(\cF_-)=\kappa_-(H)<\infty$. 
Factoring out $\cF_0$ one can assume, 
without loss of generality, that $\cF_0=0$. Then, let $\cK_+$ denote the completion 
of $(\cF_+;[\cdot,\cdot]_H)$ to a Hilbert space and then
$\cK=\cF_-[+]\cK_+$ is a Pontryagin space with $\kappa_-(\cK)=\kappa_-(H)$. For 
arbitrary $x\in X$, the linear 
operator $V(x)\colon \cH\ra\cK$ is defined by assigning to each vector $h\in\cH$ the
$\cH$-vector field $f\colon X\ra \cH$ given by $f(x)=h$ and $f(y)=0$ for all $y\in X$, 
$y\neq x$. Then $(\cK;V)$ is a minimal lineariasation of $H$.

The uniqueness of minimal linearisations of $H$, 
modulo unitary equivalence, follows 
from the fact that any dense linear subspace of a Pontryagin space with finite 
negative/positive signature contains a maximal negative/positive subspace and the 
continuity of isometric densely defined operators between Pontryagin spaces 
with the same negative/positive signature.

From Proposition~\ref{p:lik} (1), once one gets a minimal linearisation 
$(\cK;V)$ of $H$
one immediately obtains a reproducing kernel Pontryagin space $\cR$ with 
reproducing kernel $H$ as defined in \eqref{e:res} and \eqref{e:ipres}. The 
uniqueness of the reproducing kernel Pontryagin space $\cR$ follows from 
Proposition~\ref{p:lik} (2) and the uniqueness, modulo unitary equivalence, of the 
minimal linearisation of $H$.

There is a more direct but longer way of constructing the reproducing kernel \Kr\ 
space associated to a quasi positive semidefinite $\cH$-kernel $H$, by considering 
the vector space $\cR_0(\cH)$ linearly generated by $H_xh$, for $x\in X$ and 
$h\in\cH$ on which the inner product $[\cdot,\cdot]_{\cR_0}$ is defined as in 
\eqref{e:hega} in such a way that $(\cR_0(\cH);[\cdot,\cdot]_{\cR_0})$ is a 
reproducing kernel inner product space with $H$ its reproducing kernel. 
The inner product space $(\cR_0(\cH);[\cdot,\cdot]_{\cR_0})$ is 
nondegenerate and its negative signature is the same with $\kappa_-(H)$, hence 
finite. Then $\cR_0(\cH)$ is decomposable, hence $\cR_0(\cH)=\cR_-[+]\cR_{0,+}$, 
where $\dim(\cR_-)=\kappa_-(H)<\infty$ is negative definite, while $\cR_{0,+}$ 
is positive definite. It can be proven that $\cR_{0,+}$ has a completion to a Hilbert 
space $\cR_+$ inside $\cF(\cH)$ and then $\cR=\cR_-[+]\cR_+$ is the reproducing 
kernel \Kr\ space of $H$.

The next theorem points out a property of propagation to arbitrarily large domains of 
holomorphy for quasi semidefinite holomorphic kernels.

\begin{theorem}\label{t:propag} Let $K$ be a Hermitian holomorphic $\cH$-kernel on 
some region $\Omega$ and let $\Omega_0$ be a subregion of $\Omega$. If 
$\kappa_-(K_0)<\infty$, where $K_0$ is the restriction of $K$ to 
$\Omega_0\times \Omega_0$, then $\kappa_-(K)=\kappa_-(K_0)$.
\end{theorem}

The idea of proof is to use Cauchy's Theorem in order to prove the propagation 
property for positive semidefinite kernels, first for disks around the origin of the 
complex plane, then for union of regions for which the kernel is positive definite, and 
finally to use the decomposition of $K=K_+-K_-$, where $K_\pm$ are positive 
semidefinite kernels and $K_-$ has a reproducing Hilbert space of dimension 
$\kappa_-(K)$, that can be obtained from \eqref{e:dec}.

\section{Induced \Kr\ Spaces}\label{s:iks}
Let $( \cH,\langle\cdot,\cdot\rangle_\cH)$ 
be a  Hilbert  space and consider a bounded
selfadjoint operator $A$ on $ \cH$. A new inner product on
$ \cH$ is defined by
\begin{equation}\label{produs}[h,k]_A=\langle 
Ah,k\rangle_\cH,\quad h,k\in \cH.\end{equation}
In this section the properties of some  \Kr\  spaces
associated with this inner product are described.

By definition, a {\em  \Kr\  space
induced} by $A$ is a pair $( \cK,\Pi)$, where 
\begin{itemize}
\item[(ik1)] $ \cK$ is a  \Kr\  space and
$\Pi\in \cL( \cH, \cK)$ has dense range.
\item[(ik2)] $[\Pi x,\Pi y]_\cK=\langle Ax,y\rangle_\cH$, for all $ x,y\in \cH$.
\end{itemize}

There exist at least two main
constructions of  \Kr\  spaces induced by selfadjoint operators, and
they are related by certain unitary equivalences.
Two  \Kr\  spaces $( \cK_i,\Pi_i)$, $i=1,2$, induced by the same
$A$, are {\em unitary equivalent} if there exists a unitary operator
$U\in \cL( \cK_1, \cK_2)$ such that $U\Pi_1=\Pi_2$.

\begin{example}\label{ex:ka} \emph{The Induced \Kr\ Space $(\cK_A;\Pi_A)$.}
Let
$A$ denote the selfadjoint operator with respect to the Hilbert space
$( \cH,\langle\cdot,\cdot\rangle_\cH)$. Let $ \cH_-$ and $ \cH_+$ be
the spectral subspaces corresponding to the semiaxis $(-\infty,0)$ and,
respectively, $(0,+\infty)$ and the operator $A$.
Then one gets the decomposition
\begin{equation*} \cH= \cH_-\oplus\ker A\oplus \cH_+.\end{equation*}
Note that
$( \cH_-,-[\cdot,\cdot]_A)$ and $( \cH_+,[\cdot,\cdot]_A)$ are
positive definite inner product spaces and hence they can be completed
to Hilbert spaces $ \cK^-$ and, respectively,
$ \cK^+$. We can build the  \Kr\  space
$( \cK_A,[\cdot,\cdot]_A)$ by letting
\begin{equation}\label{kaa} \cK_A= \cK^-[+] \cK^+.\end{equation}
The operator $\Pi_A\in \cL( \cH, \cK_A)$ is, by definition, the
composition of the orthogonal projection of $ \cH$ onto
$ \cH\ominus\ker A$ with the embedding of $ \cH\ominus\ker A$ into $ \cK_A$.
With these definitions, it is readily verified that $( \cK_A,\Pi_A)$ is
a  \Kr\  space induced by $A$.

In order to take a closer look at the strong topology of the  \Kr\  space
$ \cK_A$, consider the seminorm
$ \cH\ni x\mapsto\||A|^{1/2}x\|=
\langle |A|^{1/2}x, |A|^{1/2}x \rangle ^{1/2}$. The
kernel of this seminorm is exactly $\ker A$ and the completion of
$ \cH\ominus\ker A$ with respect to this norm is exactly the space
$ \cK_A$. Moreover, the strong topology of $ \cK_A$ is induced by
the extension of this seminorm. The positive definite inner product
on $ \cK_A$
associated with the norm $\||A|^{1/2}\cdot\|$ is
$\langle|A|\cdot,\cdot\rangle$. Hence, if $A=S_{A}|A|$ is the
polar decomposition of $A$ and $S_{A}$ denotes the corresponding
selfadjoint partial isometry, then it follows that $S_{A}$ can be extended
by continuity to $ \cK_A$ and this extension
is exactly the fundamental symmetry
of $ \cK_A$ corresponding to $\langle|A|\cdot,\cdot\rangle$.
\end{example}

\begin{example}\label{ex:debra} \emph{The Induced \Kr\ Space 
$(\cB_A;\Pi_{\cB_A})$.} With notation as in Example~\ref{ex:ka},
consider the polar decomposition of $A=S_A|A|$ as in Example~\ref{ex:ka}. 
Define the
space $ \cB_A= \ran(|A|^{1/2})$ endowed with the positive definite
inner product $\langle\cdot,\cdot\rangle_{ \cB_A}$ given by
\begin{equation}\label{depos}\langle|A|^{1/2}x,|A|^{1/2}y\rangle_{ \cB_A}=\langle
P_{ \cH\ominus\ker A}x,y\rangle_\cH,\quad x,y\in \cH.\end{equation}
This positive definite inner product
is correctly defined and $( \cB_A,\langle\cdot,\cdot\rangle_{ \cB_A})$
is a Hilbert space. To see this, just note that the operator
$|A|^{1/2}\colon \cH\ominus\ker A\ra \cB_A$ is a Hilbert space unitary
operator.

On $ \cB_A$ one defines the inner product
$[\cdot,\cdot]_{ \cB_A}$ by
\begin{equation}\label{debin}[|A|^{1/2}x,|A|^{1/2}y]_{ \cB_A}=\langle
S_{A}x,y\rangle_\cH,\quad x,y\in \cH.\end{equation} Since the operator
$|A|^{1/2}$ and $S_{A}$ commute it follows that
\begin{equation*}[a,b]_{ \cB_A}=\langle S_{A}a,b\rangle_{ \cB_A},\quad
a,b\in \cB_A.\end{equation*}

The operator $S_{A}| \cB_A$ is a
symmetry. We define now a linear operator $\Pi_{ \cB_A}\colon \cH\ra \cB_A$ by
\begin{equation}\label{piba}\Pi_{ \cB_A}h=|A|h,\quad h\in \cH.\end{equation}
It follows that $( \cB_A,[\cdot,\cdot]_{ \cB_A})$ is
a  \Kr\  space induced by $A$.\end{example}

The  \Kr\  space induced by a selfadjoint operator
as in Example \ref{ex:debra} is a genuine operator range subspace. More
precisely, to any  \Kr\  space $ \cK$ continuously embedded in the  Hilbert  space
$ \cH$ one associates a selfadjoint operator $A\in \cL( \cH)$ in the
following way: let $\iota\colon \cK\ra \cH$ be the inclusion operator
which is supposed bounded and take $A=\iota\iota^\sharp\in \cL( \cH)$.
Clearly $A$ is selfadjoint and $( \cK,\iota^\sharp)$ is a  \Kr\  space
induced by $A$. Conversely, from Example \ref{ex:debra} it is easy to see
that $ \cB_A$ is a  \Kr\  space continuously embedded in $ \cH$. 

The connection between the induced  \Kr\  spaces $( \cK_A,\Pi_A)$ and
$( \cB_A,\Pi_{ \cB_A})$ is explained by the following
\begin{proposition}\label{arva} The induced  \Kr\  spaces $( \cK_A,\Pi_A)$ and
$( \cB_A,\Pi_{ \cB_A})$ are unitary equivalent, more precisely, the
mapping
\begin{equation*} \cK_A\ni x\mapsto |A|x\in \cB_A,\end{equation*} 
extends uniquely to a unitary operator $V\in \cL( \cK_A, \cB_A)$ such that
$V\Pi_A=\Pi_{ \cB_A}$.\end{proposition}

\section{Uniqueness of Induced \Kr\ Spaces}\label{s:uiks}

The two examples of induced  \Kr\  spaces described in 
Example~\ref{ex:ka} and Example~\ref{ex:debra} turned
out to be unitary equivalent. However, in general, not all possible  \Kr\  
spaces induced by a fixed selfadjoint operator are unitary equivalent. We denote 
by $\rho(T)$ the \emph{resolvent set} of the operator $T$.
\begin{theorem}\label{carunipas}
Let $A$ be a bounded selfadjoint operator in the Hilbert space $ \cH$. The
following statements are equivalent:

{\rm (i)} The  \Kr\  space induced by $A$ is unique, modulo  unitary equivalence.

{\rm (ii)} There exists $\epsilon>0$ such that either $(0,\epsilon)\subset\rho(A)$
or $(-\epsilon,0)\subset\rho(A)$.

{\rm (iii)} For some (equivalently, for any)  \Kr\  space $\{ \cK,\Pi\}$
induced by $A$, the range of $\Pi$ contains a maximal uniformly
definite subspace of $ \cK$.
\end{theorem}
The equivalence of (ii) and (iii) is a simple matter of spectral theory for bounded 
selfadjoint operators. The implication (iii)$\Ra$(i) comes from the fact that any
densely defined isometric operator whose domain contains a maximal uniformly 
definite subspace is bounded. So only the idea of the implication 
(i)$\Ra$(ii) is clarified.
 
Assuming that the statement (ii) does not hold, there exists a
decreasing sequence of numbers $(\mu_n)_{n\geq 1}$ 
with $\mu_n\in\sigma(A)$ and
$0<\mu_n<1$ for all $n\geq 1$, such that $\mu_n\ra0$ $(n\ra\infty)$, and there 
exists a decreasing sequence of numbers $(\nu_n)_{n\geq 1}$ with 
$-\nu_n\in\sigma(-A)$ and $0<\nu_n<1$ for all $n\geq 1$, such that 
$\nu_n\ra0$ $(n\ra\infty)$. Then, letting $\mu_0=\nu_0=1$, there exist
sequences of orthonormal 
vectors $\{e_n\}_{n\geq 1}$ and $\{f_n\}_{n\geq 1}$ such that
\begin{equation} e_n\in E((\mu_{n},\mu_{n-1}])\cH,\quad f_n\in
E([-\nu_{n-1},-\nu_{n}))\cH,\quad n\geq 1.\end{equation}
As a consequence, one also gets
\begin{equation} [Ae_i,f_j]=0,\quad i,j\geq 1.\end{equation}
Define the sequence $(\lambda_n)_{n\geq 1}$ by
\begin{equation*}\lambda_n=\max\{\sqrt{1-\mu_n^2},\sqrt{1-\nu_n^2}\}.
\end{equation*} 
Then $0<\lambda_n\leq 1$, $\lambda_n\uparrow 1$ $(n\ra\infty)$.

Consider the subspace $\cS_n$
of the \Kr\ space $\cK_A$, defined by
\begin{equation*}\cS_n=\CC e_n\dot+ \CC f_n,\quad n\geq 1,\end{equation*} 
and then define the
operators $U_n\in\cL(\cS_n)$,
\begin{equation*}U_n=\frac{1}{\sqrt{1-\lambda_n^2}}\left[\begin{array}{cc} 1 &
-\lambda_n \\ \lambda_n & -1 \end{array}\right], \quad n\geq 1.\end{equation*}
The operators $U_n$ are isometric in $\cS_n$.
Further, one defines the linear manifold $\cD_0$ in $\cK_A$ by
\begin{equation*}\cD_0=\bigcup_{k\geq 1}\cS_k\end{equation*} 
and note that the closure of
$\cD_0=\mathop
\bigvee\limits_{k\geq1}\{e_k,f_k\}$ is a regular subspace in $\cK_A$.

By construction, the linear manifold
\begin{equation*}\cD=\cD_+\dot+\cD_-=\ran(\Pi)\end{equation*} 
is dense in $\cK_A$,
where $A=A_+-A_-$ is the Jordan decomposition of $A$ and
$\cD_\pm=\ol{\ran(A_\pm)}$. Also,
$\cD_0\subseteq\cD$ and the following decomposition holds
\begin{equation*}\cD=\cD_0\dot+(\cD\cap\cD_0^\perp).\end{equation*} 
Then define a linear operator $U$ in $\cK_A$, with domain $\cD_0$ and the 
same range,  by
$U|\cS_n=U_n$, $n\geq 1$, and
$U|(\cD\cap\cD_0^\perp)=I|(\cD\cap\cD_0^\perp)$. The operator
$U$ is isometric, it has dense range as well as dense domain, and it is
unbounded since it maps uniformly definite subspaces into subspaces that are 
not uniformly definite.

$(\cK_A,\Pi)$ is a \Kr\ space induced by $A$.
Indeed, $\Pi\cH=U\Pi_A\cH\supseteq\cD$, and the latter
is dense in $\cK_A$. Further,
\begin{equation*}[\Pi x,\Pi y]=[U\Pi_A x,U\Pi_A y]=[\Pi_A x,\Pi_A y]=[Ax,y],\quad
x,y\in\cH.\end{equation*} Since $\Pi_A$ is bounded it follows that $\Pi$ is closed.
By the Closed Graph Principle, it follows that $\Pi $ is bounded. 
Since $U$ is unbounded it follows that
$(\cK_A,\Pi_A)$ is not unitary equivalent with $(\cK_A, \Pi)$. Thus, a contradiction 
with the assertion (i) is obtained.

Let $ \cK_k$ be two  \Kr\  spaces continuously embedded into the  \Kr\  space
$ \cH$ and denote by $\iota_k\colon \cK_k\ra \cH$ the corresponding
embedding operators, that is $\iota_kh=h$, $h\in \cK_k$, $k=1,2$. We
say that the  \Kr\  spaces $ \cK_1$ and $ \cK_2$ correspond to the same
selfadjoint operator $A$ in $ \cH$ if
$\iota_1\iota_1^\sharp=\iota_2\iota_2^\sharp=A$. If $A$ is nonnegative,
equivalently, the  \Kr\  spaces $ \cK_k$ are actually Hilbert spaces, this
implies that $ \cK_1= \cK_2$.  As a
consequence of Theorem \ref{carunipas} the following corollary is obtained.

\begin{corollary}\label{hara} Given a selfadjoint operator $A\in \cL( \cH)$, the
following statements are mutually equivalent:

{\rm (a)} There is a unique  \Kr\  space $ \cK$ continuously embedded in
$ \cH$ and associated to $A$.

{\rm (b)} There exists $\epsilon>0$ such that either
$(-\epsilon,0)\subset\rho(A)$ or $(0,\epsilon)\subset\rho(A)$.

{\rm (c)} There exists a  \Kr\  space $ \cK$ continuously embedded in
$ \cH$, $\iota\colon \cK\ra \cH$ such that $\iota\iota^\sharp=A$ and
$ \ran(\iota^\sharp)$ contains a maximal uniformly definite subspace of
$ \cK$.\end{corollary}

Indeed, let $\cK_i$ be two \Kr\ spaces continuously embedded in
$\cH$ and let $\iota_i\colon\cK\ra\cH$ be the embedding operators,
$i=1,2$. Assume that the induced \Kr\ spaces $(\cK_i,\iota_i^\sharp)$
are unitary equivalent, that is, there exists a unitary operator
$U\in\cL(\cK_1,\cK_2)$ such that $U\iota_2=\iota_1$. Then
$\iota_1=\iota_2U^\sharp$ and taking into account that $\iota_i$ are
embeddings, it follows that
\begin{equation*}U^\sharp x=\iota_2U^\sharp x=\iota_1x=x,\quad x\in\cK_2,
\end{equation*} hence
$\cK_1=\cK_2$ and $\iota_1=\iota_2$. This shows that the two \Kr\
spaces coincide. We can now apply Theorem \ref{carunipas} and get the
equivalence of the statements (a), (b), and (c).

\section{Existence of Reproducing Kernel \Kr\ Spaces}\label{s:erkks}

The next theorem clarifies the problem of existence of reproducing kernel \Kr\
spaces associated to Hermitian $\bah$-kernels. Notation is as in Section~\ref{s:hk}.
In addition, two positive semidefinite $\bah$-kernels are 
called \emph{independent} if
for any $P\in{\mathfrak K}^+(\bah)$ such that $P\leq H,K$, it follows $P=0$.

\begin{theorem}\label{kolmo} 
Let $H\in{\mathfrak K}^h(\bah)$. The following assertions are equivalent:
\begin{itemize}
\item[(1)] There exists $L\in{\mathfrak K}^+(\bah)$ such that $-L\leq H\leq L$.
\item[(1)$^\prime$]There exists $L\in{\mathfrak K}^+(\bah)$ such that
\begin{equation*}|[f,g]_H|\leq [f,f]_L^{1/2}\, [g,g]_L^{1/2},\quad f,g\in\cF_0(\bah).
\end{equation*}
\item[(2)] $H=H_1-H_2$ with $H_1,H_2\in{\mathfrak K}^+(\bah)$.
\item[(2)$^\prime$] $H=H_+-H_-$ with $H_\pm\in{\mathfrak K}^+(\bah)$ independent.
\item[(3)] There exists a Kolmogorov decomposition $(V;\cK)$ of $H$.
\item[(4)] There exists a \Kr\ space with reproducing kernel $H$.
\end{itemize}
\end{theorem}

Note that letting $f=g$ in (1)$^\prime$ one obtains (1). Conversely, 
let $L\in{\mathfrak K}^+(\bah)$ be such that
$-L\leq H\leq L$, that is $|[f,f]_H|\leq [f,f]_L$, for all $f\in\cF_0(\bah)$. Let
$f,g\in\cF_0(\bah)$. Since $H$ is Hermitian, one gets
$4\Re[f,g]_H=[f+g,f+g]_H-[f-g,f-g]_H$ and hence
$4\,|\Re[f,g]_H|\leq[f+g,f+g]_L+[f-g,f-g]_L=2[f,f]_L+2[g,g]_L$.
Let $\lambda\in\CC$ be chosen such that $|\lambda|=1$ and
$\Re[f,\lambda g]=[f,\lambda g]$. Then
\begin{equation}\label{felam}|[f,\lambda
g]_H|\leq\frac{1}{2}[f,f]_L+\frac{1}{2}[g,g]_L.\end{equation}
We distinguish two possible cases. First, assume that either $[f,f]_L=0$
or $[g,g]_L=0$. To make a choice assume $[f,f]_L=0$. Consider the
inequality \eqref{felam} with $g$ replaced by $tg$ for $t>0$. Then
$|[f,g]_H|\leq\frac{t}{2}[g,g]_H$. Letting $t\ra0$ one gets
$[f,g]_H=0$. 

Second case, assuming that both $[f,f]_L$ and $[g,g]_L$ are nontrivial,
in \eqref{felam} replace $f$ by $[f,f]_L^{-1/2}f$ and $g$ by
$[g,g]_L^{-1/2}g$ to get $|[f,g]_H|\leq[f,f]_L^{1/2}[g,g]_L^{1/2}$. Thus, the assertions
(1) and (1)$^\prime$ are equivalent. \smallskip

In order to describe the idea of the proof of the implication (1)$^\prime\Rightarrow$(2)$^\prime$, 
let $\cK_L$ be the quotient-completion of
$(\cF_0(\bah),[\cdot,\cdot]_L)$ to a Hilbert space.
More precisely, letting $\cN_L=\{f\in\cF_0\mid [f,f]_L=0\}$ denote
the isotropic subspace of the positive semidefinite inner product space
$(\cF_0(\bah),[\cdot,\cdot]_L)$, one considers the quotient
$\cF_0(\bah)/\cN_L$ and complete it to a Hilbert space $\cK_L$.
The inequality (1)$^\prime$ implies that the isotropic subspace $\cN_L$ is
contained into the isotropic subspace $\cN_H$ of the inner product
$(\cF_0,[\cdot,\cdot]_H)$. Therefore, $[\cdot,\cdot]_H$ uniquely
induces an inner product on $\cK_L$, also denoted by $[\cdot,\cdot]_L$
such that the inequality in (1)$^\prime$ still holds for all $f,g\in\cK_L$. By
the Riesz Representation Theorem one gets a selfadjoint and contractive operator
$A\in\cL(\cK_L)$, such that
\begin{equation}\label{fegeha} [f,g]_H=[Af,g]_L,\quad f,g\in\cK_L.\end{equation} Let
$A=A_+-A_-$ be the Jordan decomposition of $A$ in $\cK_L$. Then
$A_\pm$ are also contractions and hence
\begin{equation}\label{afefe} [A_\pm f,f]_L\leq[f,f]_L,\quad f\in\cK_L.\end{equation}
It can be proven that the nonnegative inner products $[A_\pm\cdot,\cdot]$ 
uniquely induce kernels $H_\pm\in{\mathfrak K}^+(\bah)$ such that
$[f,f]_{H_\pm}\leq[f,f]_L,\quad f\in\cF_0(\bah)$, and $H=H_+-H_-$.

Indeed, the inner product $[A_+\cdot,\cdot]$ restricted to
$\cF_0(\bah)/\cN_H$ can be extended to an inner product
$[\cdot,\cdot]_+$ on $\cF_0(\bah)$ by letting it be null onto $\cN_L$ and
hence
\begin{equation}\label{fegeplus} [f,f]_+\leq[f,f]_L,\quad f\in\cF_0(\bah).\end{equation}

Let $x,y\in X $ be arbitrary and $x\neq y$. Clearly, one can identify
the \Kr\  space $\cH_x[+]\cH_y$ with the subspace of all $\bah$-fields 
$f\in\cF_0(\bah)$ such that $\supp f\subseteq\{i,j\}$. With this identification, one
considers the restrictions of the inner products $[\cdot,\cdot]_+$ and
$[\cdot,\cdot]_L$ to $\cH_x[+]\cH_y$. The inner product
$[\cdot,\cdot]_L$ is jointly continuous with respect to the strong
topology of $\cH_x[+]\cH_y$. By (\ref{fegeplus}) and the equivalence
of (1) and (1)$^\prime$ one concludes that  the inner product $[\cdot,\cdot]_+$
is also jointly continuous with respect to the strong topology of
$\cH_x[+]\cH_y$ and hence, by the Riesz Representation Theorem,
there exists a selfadjoint operator $S\in\cL(\cH_x[+]\cH_y)$ such
that
$$[f,g]_+=[Sf,g]_{\cH_x[+]\cH_y},\quad f,g\in\cH_x[+]\cH_y.$$
Define $H_+(x,y)=P_{\cH_x}S|\cH_y$,
$H_+(x,x)=P_{\cH_x}S|\cH_x$,
$H_+(j,j)=P_{\cH_y}S|\cH_y$ and
$H_+(j,i)=P_{\cH_y}S|\cH_x=H_+(x,y)^\sharp$.

In this way one obtains a kernel $H_+\in{\mathfrak K}^h(\bah)$ such that
$H_+\leq L$ and $$[f,g]_{H_+}=[f,g]_+,\quad f,g\in\cF_0(\bah).$$ Since
the inner product $[\cdot,\cdot]_+$ is
nonnegative it follows that $H_+\in{\mathfrak K}^+(\bah)$.

Similarly one constructs the kernel $H_-\in{\mathfrak K}^+(\bah)$ such that
$H_-\leq L$ and $$[f,g]_{H_-}=[f,g]_-,\quad f,g\in\cF_0(\bah),$$ where
the inner product $[f,g]_-$ is the extension of the restriction of the
inner product $[A_-f,g]$ to $\cF_0(\bah)/\cN_L$, by letting it be
null onto $\cN_H$.

From $A=A_+-A_-$, (\ref{fegeha}) and the constructions of the kernels
$H_+$ and $H_-$ one concludes that $H=H_+-H_-$.

Let $P\in{\mathfrak K}^+(\bah)$ be such that $P\leq H_\pm$. Then
\begin{equation}\label{fefepe} [f,f]_P\leq [f,f]_L,\quad f\in\cF_0(\bah).\end{equation} 
As before,
$[\cdot,\cdot]_P$ induces a nonnegative inner product $[\cdot,\cdot]_P$
on $\cK_L$ such that (\ref{afefe}) holds for all $f\in\cK_L$. From
$P\leq H_\pm$ one concludes that
$$[f,f]_P\leq[A_\pm f,f]_L,\quad f\in\cK_L,$$ and, since $A_+ A_-=0$
this implies $[f,f]_P=0$ for all $f\in\cK_L$. Since by (\ref{fefepe})
one gets $\cN_L\subseteq\cN_P$ this implies that the inner product
$[\cdot,\cdot]_P$ is null onto the whole $\cF_0(\bah)$ and hence
$P=0$.

The implications (2)$^\prime\Rightarrow$(2) and (2)$\Rightarrow$(1) are clear. 

The most interesting implication is (1)$^\prime\Rightarrow$(3). 
In a fashion similar to (1)$^\prime\Rightarrow$(2)$^\prime$, one considers the
quotient-completion Hilbert space $\cK_L$, the representation
(\ref{fegeha}) and the Jordan decomposition $A=A_+-A_-$. The latter
yields in a canonical way a \Kr\  space $(\cK,[\cdot,\cdot]_H)$. We again
consider $\cN_L$ and $\cN_H$, the isotropic spaces
of the inner product spaces $(\cF_0(\bah),[\cdot,\cdot]_L)$ and,
respectively, $(\cF_0(\bah),[\cdot,\cdot]_H)$. From the inequality (1)$^\prime$
one gets $\cN_L\subseteq\cN_H$.

For every $x\in X $ and every vector
$h\in\cH_x$ one considers the function $h\in\cF_0(\bah)$ defined by
\begin{equation}\label{hijel}h(y)=\left\{\begin{array}{ccc} h, & & y=x, \\ 0, & & y\neq
x.\end{array}\right.\end{equation} This identification of vectors with functions in
$\cF_0(\bah)$ yields a natural embedding
$\cH_x\hookrightarrow\cF_0(\bah)$. With this embedding one defines
linear operators $V(x)\colon\cH_x\ra\cK$ by
\begin{equation*}V(x)h=h+\cN_H\in\cF_0(\bah)/\cN_H\subseteq\cK,\quad h\in\cH_x.
\end{equation*}
It follows that the linear operators $V(x)$ are bounded, for all $x\in X$, and that
$(\cK;V)$ is a minimal linearisation of the $\bah$-kernel H.

The idea of the proof of (3)$\Rightarrow$(1) deserves an explanation as well. 
Let $(\cK,[\cdot,\cdot])$ be a \Kr\  space and
$\{V(x)\}_{x\in X }$ be a family of bounded linear operators
$V(x)\in\cL(\cH_x,\cK)$, $x\in X $, such that
\begin{equation*}H(x,y)=V(y)^\sharp V(x),\quad x,y\in X .\end{equation*} 
Fix on $\cK$ a
fundamental symmetry $J$ and for each $x\in X $ fix a fundamental
symmetry $J_x$ on $\cH_x$. Then, defining the kernel $L$ by
\begin{equation*}L(x,y)=J_yV(y)^*V(x),\quad x,y\in X ,\end{equation*} 
it can be proven that $L\in{\mathfrak K}^+(\bah)$ and that $-L\leq H\leq L$.

\begin{example}\label{ex:contra} (1) Let $\cE$ be a reflexive real Banach space 
which is not a Hilbert space. For example, one can take $1<p<2$ and $\cE=
\ell^p_\RR$. Let $\cE^\prime$ denote its topological dual space. On $X
=\cE\times\cE^\prime$ consider the Hermitian form
\begin{equation}\label{e:has} H((e,\phi),(f,\psi))=\phi(f)+\psi(e),\quad e,f\in\cE,\ 
\phi,\psi\in\cE^\prime.
\end{equation}
The Hermitian form $H$
can be viewed as a Hermitian scalar kernel on $X\times X$ and it can be proven 
that it cannot be written as a difference of two positive semidefinite scalar kernels. 
Briefly, the idea is that $X$ is a selfdual Banach space when given the norm 
$\|(e,\psi)\|_X^2=\|e\|_\cE^2+\|\psi\|_{\cE^\prime}^2$, and $H$ is jointly 
continuous with 
respect to this norm, but the topological inner product space $(X;H;\|\cdot\|_X)$ is
not decomposable. 
\end{example}

\section{Uniqueness of Reproducing Kernel  \Kr\ Spaces}\label{s:urkks}

Quasi semidefinite Hermitian kernels are associated with 
reproducing kernel \Kr\ spaces which are unique, equivalently, they have 
linearisations having the uniqueness property modulo unitary equivalence.
In the general Hermitian case, the existence of a reproducing kernel \Kr\ space does 
not imply that it is unique, equivalently, the existence of a linearisation does not imply 
its uniqueness, modulo unitary equivalence. Recall that,
two linearisations $(\cK;V)$ and
$(\cH,U)$ of the same $\bah$-kernel $H$ are {\em unitary
equivalent} if there exists a unitary operator $\Phi\in\cL(\cK,\cH)$
such that for all $x\in X$ one gets $U(x)=\Phi V(x)$.

Let $H$ be an $\bah$-kernel. If $L\in{\mathfrak K}^+(\bah)$
is such that $-L\leq H\leq L$ then
one denotes by $\cK_L$ the
quotient completion of $(\cF_0(\bah),[\cdot,\cdot]_L)$ to a Hilbert
space and by $A=A_L\in\cL(\cK_L)$ the \emph{Gram operator} of the inner
product $[\cdot,\cdot]_H$ with respect to the positive semidefinite inner product
$[\cdot,\cdot]_L$, that is, $[h,k]_H=[A_Lh,k]_L$ for all $h,k\in\cK_L$.

\begin{theorem}\label{unicita} Let $H$ be an $\bah$-kernel on a set $X$
which has a minimal linearisation, equivalently, it is associated to a
reproducing kernel of a reproducing kernel \Kr\ 
space on $X$. The following assertions are equivalent:\smallskip

{\rm (i)} The $\bah$-kernel $H$ has unique minimal linearisation,
modulo unitary equivalence.\smallskip

{\rm (ii)} For any (equivalently, there exists a) positive semidefinite
$\bah$-kernel $L$ such that $-L\leq H\leq L$ there exists $\epsilon>0$ such that
either $(0,\epsilon)\subset\rho(A_L)$ or $(-\epsilon,0)\subset\rho(A_L)$.

{\rm (iii)}  $H$ has a minimal linearisation (equivalently, any minimal linearisation)
$(\cK;V)$ that has
fundamental decomposition $\cK=\cK^+[+]\cK^-$ such that either
$\cK^+$ or $\cK^-$ is contained in the linear manifold generated by
$V(x)\cH_x$, $x\in X$.

{\rm (iv)} The reproducing kernel \Kr\ space with reproducing kernel $H$ is unique.
\end{theorem}

In order to explain the implication (i)$\Ra$(ii), assume that there exists a positive 
semidefinite
$\bah$-kernel $L$ such that $-L\leq H\leq L$ and for any $\epsilon>0$ one gets
$(0,\epsilon)\cap\sigma(A_L)\neq\emptyset$ and
$(-\epsilon,0)\cap\sigma(A_L)\neq\emptyset$.
From Theorem \ref{carunipas} it follows that there exists two \Kr\ spaces
$(\cK,\Pi)$ and $(\cH,\Phi)$ induced by the same selfadjoint operator
$A_L$, which are not unitary equivalent. It is easy to see that the
operator $\Psi\colon\ran(\Pi)\ra\ran(\Phi)$ defined by
$$\Psi\Pi f=\Phi f,\quad f\in\cK_L,$$
is isometric, densely defined with dense range, and it is unbounded
due to the non-unitary equivalence of the two induced \Kr\ spaces.
As a consequence, $\Psi$ is closable and its closure, denoted
also by $\Psi$, shares the same properties.

Let $(V;\cK)$ be the minimal linearisation of $H$ defined as in the
proof of Theorem \ref{kolmo}, (1)$^\prime\Ra$(3).
Define a new minimal linearisation
$(\cH;U)$ and prove that it is not unitary equivalent with
$(V;\cH)$. More precisely, let $U(x)=\Psi V(x)$ for all $x\in X$. Since
$\ran(V(x))\subseteq\cD(\Psi)$ and $\Psi$ is closed it follows, via the
Closed Graph Principle, that $U(x)\in\cL(\cH_x,\cK)$ for all
$x\in X$.

Let $x,y\in X$ be arbitrary and fix vectors $h\in\cH_x$ and
$k\in\cH_y$. Then
\begin{equation*}[U(y)k,U(x)h]=[\Psi V(y)k,\Psi V(x)h]=[V(y)k,V(x)h]=[H(x,y)k,h],
\end{equation*} and hence
${U}_x^\sharp {U}_y=H(x,y)$. Also,
\begin{equation*}\bigvee\limits_{j\in X}U(y)\cH_y=\bigvee\limits_{j\in X}\Psi
V(y)\cH_y=\cl\ran(\Phi)=\cH.\end{equation*}
Thus, $(\cH;U)$ is a minimal linearisation of
the $\bah$-kernel $H$. On the other hand, since the operator $\Psi$
is unbounded it follows that the two minimal linearisations
$(\cK;V)$ and $(\cH;U)$ are not unitary equivalent.

For the implication (ii)$\Ra$(i), let $(\cK;V)$ and
$(\cH;U)$ be two minimal linearisations of $H$.
Let $J$ and $J_x$ be fundamental symmetries on $\cK$ and,
respectively, $\cH_x$, $x\in X$. We consider the positive definite
$\bah$-kernel
$L_V$ defined by $$L_V(x,y)=J_yV(y)^*V(x),\quad x,y\in X,$$ and as in
the proof of Theorem \ref{kolmo} it follows that $-L_V\leq H\leq L_V$.
We define a linear operator $\Pi_V\colon\cF_0(\bah)\ra\cK$ by
\begin{equation}\label{pive}\Pi_V(h)=\sum_{x\in X}V(x)h_x,\quad
h=(h_x)_{x\in X}\in\cF_0(\bah).\end{equation} Taking into account of the
axiom (kd2) in the definition of a minimal linearisation one obtains
\begin{equation}\label{harpv}[\Pi_Vh,\Pi_Vk]_H=[h,k]_\cK,\quad h,k\in\cF_0(\bah),\end{equation}
that is, the operator $\Pi_V$ is isometric
from $(\cF_0(\bah),[\cdot,\cdot]_H)$ into $(\cK,[\cdot,\cdot]_\cK)$.
In addition, $\Pi_V$ is also isometric when considered as
a linear operator from
$(\cF_0,[\cdot,\cdot]_{L_V})$.

Similarly, considering the positive semidefinite $\bah$-kernel $L_U$ defined by
\begin{equation*}L_U(x,y)=J_yU(y)^*U(x),\quad x,y\in X,\end{equation*} 
one gets $-L_U\leq H\leq
L_U$ and, defining the linear operator $\Pi_U\colon\cF_0(\bah)\ra\cH$ by
\begin{equation}\label{piu}\Pi_U(h)=\sum_{x\in X}U(x)h_x,\quad
h=(h_x)_{x\in X}\in\cF_0(\bah),\end{equation} one obtains
\begin{equation}\label{harpu}[\Pi_Uh,\Pi_Uk]_H=[h,k]_\cK,\quad h,k\in\cF_0(\bah),\end{equation}
that is, the operator $\Pi_U$ is isometric
from $(\cF_0,[\cdot,\cdot]_H)$ into $(\cK,[\cdot,\cdot]_\cH)$ and
$\Pi_U$ is also isometric when considered as
a linear operator from
$(\cF_0,[\cdot,\cdot]_{L_U})$ into $(\cK,\langle\cdot,\cdot\rangle_J)$.

Let $L=L_V+L_U$ and clearly $-L\leq H\leq L$. Since
$L_V\leq L$ it follows that $\cK_L$ is contractively embedded into
$\cK_{L_V}$ and hence $\Pi_V$ induces a bounded operator
$\Pi_V\colon\cA_L\ra\cK$. From (\ref{harpv}) it follows
\begin{equation} \label{kalv} \Pi_V^*J \Pi_V=A_L.\end{equation}
By assumption, there exists $\epsilon>0$
such that either $(-\epsilon,0)\subset\rho(A_L)$ or
$(0,\epsilon)\subset\rho(A_L)$ and hence from (\ref{kalv})
and taking into account that by the minimality axiom (b)
of the minimal linearisation the operator $\Pi_V$ has dense range,
it follows that there exists a uniquely determined
unitary operator $\Phi_V\colon\cK_{A_L}\ra\cK$ such that
\begin{equation}\label{philv} \Phi_Vh=\Pi_Vh,\quad h\in\cA_L,\end{equation} where $\cK_{A_L}$
is the \Kr\ space induced by the operator $A_L$.

Similarly, performing the same operations with respect to the linearisation 
$(\cH;V)$ one gets a uniquely determined
unitary operator $\Phi_U\colon\cK_{A_L}\ra\cH$ such that
\begin{equation}\label{philu} \Phi_Uh=\Pi_Uh,\quad h\in\cA_L.\end{equation}
Define the unitary operator $\Phi\colon\cK\ra\cH$ by
\begin{equation*}\Phi=\Phi_U\Phi_V^{-1}.\end{equation*} 
Taking into account of (\ref{philv}) and
(\ref{philu}), the definition of the operator $\Pi_V$ as in
(\ref{pive}), and the definition of $\Pi_U$ as in (\ref{piu}), it follows
that
\begin{equation*}\Phi(\sum_{x\in X}V(x)h_x)=\sum_{x\in X}U(x)h_x,\quad
(h_x)_{x\in X}\in\cF_0(\bah).\end{equation*} 
This implies readily that for all
$x\in X$ one has $\Phi V(x)=U(x)$ and hence the two Kolmogorov
decompositions $(\cK;V)$ and
$(\cH;U)$ are unitary equivalent.\medskip

As a consequence of Theorem \ref{unicita} one can obtain a rather general 
sufficient condition of nonuniqueness. Let
$K$ and $H$ be two positive semidefinite $\bah$-kernel. Then
one considers the Hilbert spaces $\cK_K$ and $\cK_H$, obtained by
quotient completion of $(\cF_0,[\cdot,\cdot]_K)$ and, respectively, of
$(\cF_0,[\cdot,\cdot]_H)$. If $H\geq K$ then $\cK_H$ is contractively
embedded into $\cK_K$. The kernel $H$
is \emph{$K$-compact} if the embedding of $\cK_H$ into $\cK_K$ is a compact
operator.
\begin{corollary}\label{schwartz} Let $H_+,H_-\in\cK^+(\bah)$ be two independent
kernels, both of them of infinite rank. If there exists a kernel
$K\in\cK^+(\bah)$ such that $H_+$ and $H_-$ are $K$-compact, then the
minimal linearisations of the kernel $H_+-H_-$ are not unique,
modulo unitary equivalence.\end{corollary}

 Let $H=H_1-H_2$. Clearly $-K\leq H\leq K$. Let $A_\pm\in\cL(\cK_K)$
denote the Gram operator of the kernel $H_\pm$. Since $H_+$ and $H_-$
are independent it follows that $A=A_+-A_-$ is the Gram operator of $H$.
Since $H_\pm$ are of infinite rank and $K$-compact it follows that
$A_\pm$ are compact operators of infinite rank in $\cL(\cK_K)$ and
hence the spectra $\sigma(A_\pm)$ are accumulating to $0$. Then the
spectrum $\sigma(A)$ is accumulating to $0$ from both sides. This
clearly contradicts the condition (ii) in Theorem \ref{unicita} and hence
the kernel $H$ has non-unique minimal linearisations.\qed\medskip

\section{Holomorphic Kernels: Single Variable Domains}\label{s:hksvd}

Let $\Omega$ be a domain, a nonempty open subset, in the complex field $\CC$, 
and let $\cH$ be a Hilbert space. 
A kernel $K\colon \Omega\times \Omega\ra \cL(\cH)$ 
is called \emph{holomorphic} if it is holomorphic in the first variable, that is, for each 
$w\in \Omega$, the map $\Omega\ni z\mapsto K(z,w)\in\cL(\cH)$ and 
conjugate holomorphic in the second variable, that is, for each $z\in \Omega$, 
the map $\ol{\Omega}\ni \ol{w}\mapsto K(z,\ol{w})\in \cL(\cH)$ is holomorphic. Recall
that, for Banach space valued functions of complex variable, strong holomorphy 
is the same with weak holomorphy. If $K$ is Hermitian, then $K$ is a 
holomorphic 
kernel if and only if the map $\Omega\ni z\mapsto K(z,w)\in\cL(\cH)$ 
is holomorphic for all $w\in \Omega$.

\begin{theorem}\label{t:alpay} Let  
$\cH$ be a Hilbert space and, for some $r>0$, let $K$ be a holomorphic 
$\cH$-kernel on $\DD_r=\{z\in\CC\mid |z|<r\}$. Then, there exists $0<r^\prime\leq r$ 
and a reproducing kernel \Kr\ space on $\DD_{r^\prime}$ with reproducing kernel 
$K|\DD_{r^\prime}\times \DD_{r^\prime}$.
\end{theorem}

The first step in the proof of this theorem is to observe that, without loss of 
generality, one can assume $r>1$. Indeed, if $r\leq 1$ then, for some 
$0<\rho<r$ small enough, the $\cH$-kernel $K_\rho(z,w)=K(\rho z,\rho w)$ 
is holomorphic on $\DD_{r/\rho}$. If $\cK_\rho$ is the reproducing kernel \Kr\ space
with reproducing kernel $K_\rho$ restricted to $\DD_{r^{\prime\prime}}$, for some 
$0<r^{\prime\prime}\leq r/\rho$, let $r^\prime=r^{\prime\prime}\rho$ and let 
$\cK$ denote the vector space of functions $f\colon \DD_{r^\prime}\ra \cH$ such that
$f(z)=F(z/\rho)$ for some $F\in \cK_\rho$ and all $z\in \DD_{r^\prime}$. On $\cK$ 
there is defined the inner product $[\cdot,\cdot]_\cK$
\begin{equation*}[f,g]_\cK=[F,G]_{\cK_\rho},\quad f(z)=F(z/\rho),\ g(z)=G(z/\rho),\ 
F,G\in \cK_\rho.
\end{equation*} Then $(\cK;[\cdot,\cdot]_\cK)$ is a \Kr\ space. For each 
$w\in\DD_{r^\prime}$ the map $z\mapsto K(z,\rho w)$ belongs to $\cK$ and, 
for each $F\in \cK_\rho$, $f(z)=F(z/\rho)$, and $h\in\cH$,
\begin{equation*} [f,K(\cdot,\rho w)h]_\cK=[F,K(\rho\cdot,\rho w)h]_{\cK_\rho}=
[F(\rho w),h]_\cH=[f(w),h]_\cH,
\end{equation*} hence $\cK$ is a reproducing kernel \Kr\ space with reproducing 
kernel $K|\DD_{r^\prime}\times\DD_{r^\prime}$.

There are two main ideas of the proof. First, 
the Szeg\"o kernel $S$, see Example~\ref{ex:hardy}, plays a distinguished role in 
holomorphy, and allows us to construct the convolution kernel of $K$ in the Hardy 
space $H^2(\DD)$. 
The second idea is that, once the convolution operator represented 
as a selfadjoint bounded operator on a Hilbert space of functions is defined, 
the construction of 
the induced 
\Kr\ space as in Example~\ref{ex:debra} will provide the reproducing 
kernel \Kr\ space with reproducing kernel $K$. Here are a few details.

Letting $r>1$, there is considered the Hardy space $H^2(\DD)\otimes \cH$, 
identified with a space of $\cH$-valued functions $f(z)=\sum_{n=0}^\infty a_n z^n$,
where $(a_n)_{n\geq 0}$ is a sequence of vectors in $\cH$ such that 
$\|f\|^2=\sum_{n\geq 0}^\infty  \|a_n\|^2_\cH<\infty$.  Also,
the inner product on the Hilbert space $H^2(\DD)\otimes \cH$ is
\begin{equation*} \langle f,g\rangle=\frac{1}{2\pi} \int_0^{2\pi} 
\langle f(\emath^{\iac t}), g(\emath^{\iac t})\rangle_\cH\de t, 
\quad f,g\in H^2(\DD)\otimes \cH.
\end{equation*}
Thus, on $H^2(\DD)\otimes \cH$
one can define the analog of the convolution operator $K$
\begin{equation}\label{e:ped} (Kf)(z)=\frac{1}{2\pi}\int_0^{2\pi} 
K(z,\emath^{\iac t})f(\emath^{\iac t}) \de t,\quad f\in H^2(\DD)\otimes \cH.
\end{equation}
Letting $M=\sup_{|z|,|w|\leq 1} \|K(z,w)\|<\infty$, it follows that 
$\|Kf\|\leq M \|f\|_{H^2(\DD)\otimes \cH}$, hence the convolution operator $K$ is 
a bounded linear operator in $H^2(\DD)\otimes \cH$. On the other hand,
\begin{equation*}\langle Kf,g\rangle_{H^2(\DD)\otimes \cH}
=\frac{1}{4\pi^2} \int_0^{2\pi} \int_0^{2\pi} \langle K(\emath^{\iac t},\emath^{\iac s})
f(\emath^{\iac s}),g(\emath^{\iac t})\rangle_\cH\de t\de s,
\end{equation*} hence $K$ is selfadjoint. Since, for any $w\in \DD_1$ and any 
$h\in \cH$, the function $f_w(z)=h/(1-\ol{w}z)$ belongs to $H^2(\DD)\otimes \cH$ 
and by the Cauchy formula, $(Kf_w)(z)=K(z,w)h$, the range of the 
convolution operator $K$ contains all the functions $K_w(\cdot)h$, for $w\in\DD$ 
and $h\in\cH$. 

Then, one can use the construction of the induced \Kr\ space $(\cB_K;\Pi_{\cB_K})$
inside of $H^2(\DD)\otimes \cH$,
as in Example~\ref{ex:debra} but applied to the 
bounded selfadjoint operator $K$ in the Hilbert space $H^2(\DD)\otimes \cH$, 
in order to get the reproducing kernel \Kr\ space $\cK$ with reproducing
kernel $K$. 

\section{Holomorphic Kernels: Several Variables Domains}\label{s:hksevd}

In this section it is considered the analog of Theorem~\ref{t:alpay} in case the
Hermitian kernels are defined on domains in $\CC^N$ for $N\geq 2$. We use the 
notation as in Example~\ref{ex:das} where the Drury-Arveson space was 
constructed as the reproducing kernel Hilbert space associated to the Szeg\"o 
kernel. For simplicity, it is considered only scalar valued kernels.

Recall that a scalar-valued Hermitian kernel $K$, 
defined on a nonempty open subset $\cO$ of $\cG=\CC^N$,
is {\it holomorphic on} $\cO$ if $K(\cdot ,\eta )$ is holomorphic on $\cO$ for
each fixed $\eta \in \cO$. Since $K$ is Hermitian, it follows that $K$ is conjugate 
holomorphic in the second variable.

\begin{theorem}\label{t:galpay} Let $r>0$ and let $K$ be a Hermitian
holomorphic kernel on the open ball $\BB_r$ in $\CC^N$. Then, 
there exist $0<r^\prime\leq r$ and a reproducing kernel \Kr\
space $\cK$ on $\BB_{r^\prime}$ with reproducing kernel 
$K|\BB_{r^\prime}\times\BB_{r^\prime}$.
\end{theorem}

To a certain extent, the proof of this theorem follows a pattern similar to that 
of the proof of Theorem~\ref{t:alpay}, namely, first a scaling argument can be used
in order to reduce the proof to the case $r>1$, then
the convolution kernel of $K$ can defined on the Drury-Arveson space 
$H^2(\BB_1)$ associated to the Szeg\"o kernel $S$, see Example~\ref{ex:das}, and
it can be proven that this convolution operator is bounded and selfadjoint. Finally,
the construction of type $(\cB_A;\Pi_{\cB_A})$, see Example~\ref{ex:debra}, 
can be used in order to produce
a reproducing kernel space $\cK$ with reproducing kernel 
$K$. This proof shows, once again, a certain 
universality property of the Szeg\"o kernel with respect to holomorphic Hermitian 
kernels.

Let $K$ be a scalar Hermitian holomorphic kernel on $B_r$, $r>0$. 
Since $K$ is Hermitian, it follows that $\zeta \mapsto K(\xi ,\ol\zeta )$ is 
holomorphic on $B_r$ for each $\xi \in B_r$, that is, letting $\{e_{j}\}_{j=1}^N$
denote the canonical orthonormal basis of $\cG=\CC^N$, the 
conjugation can be defined by
\begin{equation*}\xi =\sum _{j=1}^N\langle \xi ,e_{j}\rangle e_{j}
\rightarrow 
\sum_{j=1}^N\overline{\langle \xi ,e_{j }\rangle }e_{j }=
\overline\xi,\end{equation*}
so that the function $f(\xi ,\eta )=K(\xi,\ol\eta )$
is separately holomorphic on $B_r\times B_r$. By Hartogs' Theorem $f$ is 
holomorphic on $B_r\times B_r$, hence $f$ is locally bounded.
Similar to the argument provided for Theorem~\ref{t:alpay}, without loss of 
generality one can suppose that $r>1$. Hence
there exist $1<\rho <r$ and $C>0$ such that:
\begin{equation}\label{cinciunu}
|K(\xi ,\eta )|\leq C \quad \mbox{for all}\quad \xi ,\eta \in B_{\rho }
\end{equation}
and 
\begin{equation}\label{cincidoi}
K(\xi,\ol\eta )=\sum _{m\geq 0}p_m(\xi, \eta)
\end{equation}
uniformly on $B_{\rho }$, where each $p_m$, $m\geq 0$, is an 
$m$-homogeneous complex polynomial on $2N$ variables. 
There exists a continuous linear functional $A_m$ on 
$P_m(\cG\times\cG)^{\otimes m}$, see \eqref{e:peme}, such that 
\begin{equation}\label{cincitrei}
p_m(\xi ,\eta )=A_m((\xi ,\eta )^{\otimes m}), \quad\mbox{for all }\xi,\eta\in\CC^N.
\end{equation}
Using Cauchy Inequalities, for $B_{\rho }$, one gets
\begin{equation}\label{cincipatru}
\|A_m\|\leq C\displaystyle\frac{1}{\rho ^m},
\end{equation}
hence
\begin{equation}\label{cincicinci}
\sum _{m\geq 0}\|A_m\|^2\leq C\sum _{m\geq 0}
\displaystyle\frac{1}{\rho ^{2m}}=C
\displaystyle\frac{1}{1-1/\rho ^2}=C'<\infty.
\end{equation}
 By the Riesz Representation Theorem, there exist
$a_m\in P_m(\cG \times \cG )^{\otimes m}$, $m\geq 0$, such that 
\begin{equation}\label{cincisase}
A_m((\xi ,\eta )^{\otimes m})=\langle (\xi ,\eta )^{\otimes m},a_m\rangle_{(\cG\times\cG)^{\otimes m}},
\end{equation}
and 
\begin{equation}\label{cincisapte}
\|a_m\|=\|A_m\|,
\end{equation}
(with $a_0=A_0\in \CC$).
Since $P_m(\cG \times \cG )^{\otimes m}$ is isometrically isomorphic 
to $(P_m\cG ^{\otimes m})^{\oplus (m+1)}$, it is deduced that 
there are $a_m^k\in P_m\cG ^{\otimes m}$, $k=0,\ldots ,m$, such that
\begin{equation}\label{cinciopt}
\langle (\xi ,\eta )^{\otimes m},a_m\rangle_{(\cG\times\cG)^{\otimes m}}=
\sum _{k=0}^m\langle b_m^k(\xi ,\eta ),a_m^k\rangle_{\cG^{\otimes m}},
\end{equation}
and 
\begin{equation}\label{cincinoua} 
\sum _{k=0}^m\|a_m^k\|^2=\|a_m\|^2,
\end{equation}
where $b_0^0=1$ and $b_m^k(\xi,\eta)=\xi ^{\otimes (m-k)}\otimes
\eta ^{\otimes k}$, $m\geq 1$, $k=0,\ldots ,m$.
By \eqref{cincidoi}, 
\eqref{cincitrei}, \eqref{cincisase}, and \eqref{cinciopt},
\begin{equation*}K(\xi ,\eta )=\sum _{m\geq 0}\sum _{k=0}^m
\langle b_m^k(\ol\xi,\eta ),a_m^k\rangle =\sum _{k\geq 0}\sum _{m\geq k}
\langle b_m^k(\ol\xi,\eta ),a_m^k\rangle ,\end{equation*}
where the series converge absolutely on $\eta $ by \eqref{cincipatru}.

Using all these, it can be shown that $K_\eta\in H^2(\BB_1)$ for all 
$\eta\in \BB_1$, where $K_\eta(\xi)=K(\xi,\eta)$. Then letting $Ka_\eta=K_\eta$, 
$\eta\in\BB_1$, one gets a bounded linear operator in $H^2(\BB_1)$ such that
\begin{equation*}K(\xi,\eta)=\langle Ka_\xi,a_\eta\rangle_{H^2(\BB_1)}.
\end{equation*} This operator $K$ is selfadjoint and it is the analog of 
the convolution operator for which, applying the construction of type 
$(\cB_A;\Pi_{\cB_A})$ as in Example~\ref{ex:debra}, 
one gets a reproducing kernel \Kr\ space $\cK$ 
with reproducing kernel $K|\BB_1\times\BB_1$.

\section{Comments}\label{s:c}

The theory of reproducing kernel Hilbert spaces and their positive semidefinite 
kernels originates with the works of  S.~Zaremba 
\cite{Zaremba}, G.~Szeg\"o \cite{Szego}, S.~Bergman \cite{Bergman}, 
and S.~Bochner \cite{Bochner}. E.H.~Moore \cite{Moore} also 
contributed significantly to this theory, but the first systematisation and abstract 
presentation belongs 
to N.~Aronszajn \cite{Aronszajn}. A different but equivalent theory belongs to 
L.~Schwartz \cite{Schwartz}, whose work remained almost unnoticed for a long time, 
although it was the first to consider reproducing kernel \Kr\ spaces (Hermitian 
spaces, as called there). So far, monographs 
on this subject have been written by T.~Ando \cite{Ando}, S.~Saitoh \cite{Saitoh}, 
and the forthcoming title of S.~Saitoh and Y.~Sawano \cite{SaitohSawano}, which 
are good sources for the large area of applications of the technique of 
reproducing kernel spaces in complex functions theory, ordinary and partial
differential equations, integral equations, and approximation and numerical analysis.

So far, there are two monographs devoted to indefinite inner product spaces and 
their linear operators, J.~Bognar \cite{Bognar} and T.Ya.~Azizov and I.S.~Iokhvidov
\cite{AzizovIohvidov}, 
where proofs, examples, and counter-examples of the facts recalled in Section{s:ks} 
can be found.

The introductory material on Hermitian kernels in Section~\ref{s:hk} follows 
T.~Constantinescu and A.~Gheondea \cite{CG2}. Theorem~\ref{t:reka} is classical.
The concept of linearisation originates with J.~Mercer \cite{Mercer}, for the 
scalar case, and A.N.~Kolmogorov \cite{Kolmogorov}, for the operator valued case.

The Hardy Space $H^2(\DD)$ originates with the G.~Szeg\"o Kernel \cite{Szego}. 
For the theory of Hardy spaces there are monographs of P.L.~Duren 
\cite{Duren} and P.~Koosis \cite{Koosis}. Bergman kernel is also another important 
example of a positive semidefinite kernel but it falls out of this chapter concern, see
P.L.~Duren and B.~Schuster \cite{DurenSchuster}. The Drury-Arveson space 
originates with S.W.~Drury \cite{Drury} and W.B.~Arveson \cite{Arveson}.
The holomorphic kernels considered at Example~\ref{ex:sc} make the main 
object of investigation of the monograph of D.~Alpay et al.\ \cite{ADRS}. The
investigations of A.V.~Potapov \cite{Potapov}, 
L.~de~Branges \cite{dB1}--\cite{dB3}, M.G.~Kre\u\i n and 
H.~Langer \cite{KL1}--\cite{KL5}, and H.~Dym \cite{Dym} 
highly motivate the interest for Hermitian kernels
with or without finite negative signatures.
The study of T\"oplitz type kernels is related to the investigations on operator
dilations of M.A.~Na\u\i mark \cite{Naimark} and B.~Sz.-Nagy \cite{BSzNagy}. 
Our short 
presentation as in Example~\ref{ex:tk} follows \cite{CG2}. We only mention that
there is a more 
general and powerful theory of kernels invariant under actions of $*$-semigroups 
presented in \cite{CG3}, motivated by problems in mathematical physics as in 
D.E.~Evans and J.T.~Lewis \cite{EvansLewis} and K.R.~Parthasaraty, K.~Schmidt 
\cite{ParthasaratySchmidt}, and many others.

Theorem~\ref{t:qsk} essentially belongs to P.~Sorjonen \cite{Sorjonen} but
this result follows from the more general theory of L.~Schwartz \cite{Schwartz} that 
has been obtained about ten years before. Theorem~\ref{t:propag} can be found in 
D.~Alpay et al. \cite{ADRS}. 

Examples \ref{ex:ka}, \ref{ex:debra}, and Theorem~\ref{carunipas} on 
induced \Kr\ spaces can be found in \cite{CG2}, while Corollary~\ref{hara} belongs
to T.~Hara \cite{Hara}. Similar and, in a certain way, 
equivalent uniqueness conditions, 
can be found in \cite{CG1} and M.A.~Dritschel \cite{Dritschel}. B.~\'Curgus and 
H.~Langer \cite{CurgusLanger} proves that once non-equivalent induced \Kr\ 
spaces exist, there are infinitely many.

The characterisations of existence of reproducing \Kr\ spaces associated to given 
Hermitian kernels as in Theorem~\ref{kolmo} belong essentially to L.~Schwartz 
\cite{Schwartz}, but our presentation follows \cite{CG2}. Example~\ref{ex:contra} is 
from \cite{Schwartz} as well, cf.\ \cite{ADRS}.

The uniqueness Theorem~\ref{unicita} is from \cite{CG2} while 
Corollary~\ref{schwartz} is from \cite{Schwartz}. 

The result in Theorem~\ref{t:alpay} on single variable holomorphic kernels 
belongs to D.~Alpay \cite{Alpay}, while its several variables generalisation in
Theorem~\ref{t:galpay} is from \cite{CG4}: for the basics of several complex 
variables holomorphic functions facts used during the explanation of the 
ideas of the proof, see R.M. Range 
\cite{Range}.

\end{document}